\documentclass[draft]{birkau}
\usepackage{amsmath}
\usepackage{amsfonts}
\usepackage{amssymb}

\numberwithin{equation}{section}
\theoremstyle{plain}
\newtheorem{theorem}{Theorem}[section]
\newtheorem{lemma}[theorem]{Lemma}
\newtheorem{corollary}[theorem]{Corollary}
\newtheorem{proposition}[theorem]{Proposition}
\theoremstyle{definition}
\newtheorem{definition}[theorem]{Definition}
\newtheorem{example}[theorem]{Example}
\newtheorem{remark}[theorem]{Remark}
\newtheorem*{illustration}{Illustration}
\newtheorem{observation}[theorem]{Observation}
\input{tcilatex}
\begin{document}

\title[Riesz and pre-Riesz monoids]{Riesz and pre-Riesz monoids}
\author[M. Zafrullah]{Muhammad Zafrullah}
\email{mzafrullah@usa.net}

\begin{abstract}
Call a directed partially ordered cancellative divisibility mon\-oid $M$ a
Riesz monoid if for all $x,y_{1},y_{2}\geq 0$ in $M,$ $x\leq
y_{1}+y_{2}\Rightarrow x=x_{1}+x_{2}$ where $0\leq x_{i}\leq y_{i}$. We
explore the necessary and sufficient conditions under which a Riesz monoid $%
M $ with $M^{+}=\{x\geq 0|x\in M\}=M$ generates a Riesz group and indicate
some applications. We call a directed p.o. monoid $M$ $\Pi $-pre-Riesz if $%
M^{+}=M$ and for all $x_{1},x_{2},
\dots ,x_{n}\in M$, $glb(x_{1},x_{2},...,x_{n})=0$ or there is $r\in \Pi $
such that $0<r\leq x_{1},x_{2},...,x_{n},$ for some subset $\Pi $ of $M.$ We
explore examples of $\Pi $-pre-Riesz monoids of $\ast $-ideals of different
types. We show for instance that if $M$ is the monoid of nonzero (integral)
ideals of a Noetherian domain $D$ and $\Pi $ the set of invertible ideals, $%
M $ is $\Pi $-pre-Riesz if and only $D$ is a Dedekind domain. We also study
factorization in pre-Riesz monoids of a certain type and link it with
factorization theory of ideals in an integral domain.
\end{abstract}

\maketitle

$^{{}}$%
\address{Department of Mathematics, Idaho State University, Pocatello, 83209
ID}

\urladdr{https://www.lohar.com} 

\section{\textbf{Introduction\label{S1}}}

By a \emph{monoid} we shall mean a commutative semigroup with identity,
denoted 
$M=\langle M,\ast ,e\rangle $, where $\ast $ is the monoid operation and $e$
is the identity of $M$ under $\ast $. Also by a \emph{partially ordered 
\textup{(}p.o.\textup{) }monoid} we shall mean a monoid $M$ whose partial
order $\leq $ is compatible with the operation defined on the monoid, that
is $x\leq y$ implies $z\ast x\leq z\ast y$, for all $x,y,z\in M$. We denote
such monoids $M$ by $M=\langle M,\ast ,e,\leq \rangle $. A p.o. monoid $%
\langle M,\ast ,e,\leq \rangle $ is called \emph{upper \textup{(}resp., lower%
\textup{) }directed} if for each pair of elements $x,y\in M$ there is a $%
z\in M$ such that $z\geq x,y$ (resp., $z\leq x,y$). A monoid that is both
upper and lower directed is called \emph{directed}. Also $M=\langle M,\ast
,e,\leq \rangle $ is called a \emph{divisibility monoid} if $x\leq y$ in $M$
if and only if there is $a\in M^{+}=\{m\in M\mid m\geq e\}$ such that $%
y=x\ast a$. An upper directed divisibility monoid is clearly directed.
Finally we call a (p.o.) monoid \emph{cancellative} if $(z\ast x\leq z\ast y$
$\Rightarrow x\leq y)$ $z\ast x=z\ast y\Rightarrow x=y$ for all $x,y,z\in M$.

Call a directed p.o. monoid $M=\langle M,\ast ,e,\leq $ $\rangle $ a \emph{%
Riesz monoid}, if (a) $M$ is a cancellative divisibility monoid and (b)
every element $x$ of $M$ is primal i.e. for $y_{1},y_{2}\in M,$ $x\leq
y_{1}\ast y_{2}$ $\Rightarrow x=x_{1}\ast x_{2}$ such that $x_{i}\leq y_{i}$%
. It turns out that Riesz monoids have the so called \emph{Riesz
Interpolation Property}, and we indicate conditions under which a Riesz
monoid generates a \emph{Riesz group}. Since Riesz monoids have been
considered as additive as well as multiplicative monoids, under different
guises, with the same theory, we have taken $\ast $ to serve as the monoid
operation. Call elements $x_{1},...,x_{n}$ in a monoid $M=\langle M,\ast
,e,\leq $ $\rangle $ disjoint if $z\leq x_{i}$ implies $z\leq e$. We shall
use $x_{1}\wedge ...\wedge x_{n}=e$ to denote the fact that $x_{1},...,x_{n}$
are disjoint, i.e., $z\leq x_{1},...,x_{n}$ implies $z\leq e$. Call a
directed monoid $M$ a \emph{pre-Riesz monoid} if for each finite set of $%
x_{1},...,x_{n}$ in $M^{+}$, $x_{1}\wedge ...\wedge x_{n}=e$ or there is $%
t\in M$ such that $e<t\leq x_{1},...,x_{n}.$

One aim of this article is to study Riesz and pre-Riesz monoids, and how
they can be used, and the other is to study factorization in pre-Riesz
monoids and link it with the current work on factorization of ideals of an
integral domains as in \cite{AZ 2019}. Our work involves star operations on
integral domains and other tools, enough that we postpone giving the plan of
the paper until after we have given sufficient information on star
operations and related tools.

Briefly let $D$ be an integral domain with quotient field $K$, unless stated
otherwise we assume that $D\neq K$. Let $F(D)$ (resp., $f(D)$) be the set of
nonzero fractional ideals (resp., nonzero finitely generated fractional
ideals) of $D$. A star operation $\star $ on $D$ is a function on $F(D)$
that satisfies the following properties for every $I,J\in F(D)$ and $0\neq
x\in K$:


\begin{enumerate}
\item $(x)^{\star }=(x)$ and $(xI)^{\star }=xI^{\star }$,

\item $I\subseteq I^{\star }$ and $I^{\star }\subseteq J^{\star }$ whenever $%
I\subseteq J$, and

\item $(I^{\star })^{\star }=I^{\star }$.
\end{enumerate}

Now, an ideal $I\in F(D)$ is a $\star $-\emph{ideal} if $I^{\star }=I$, so a
principal ideal is a $\star $-ideal for every star operation $\star $.
Moreover $I\in F(D)$ is called a $\star $-ideal of finite type if $%
I=J^{\star }$ for some $J\in f(D)$. It can be shown that (a) for every star
operation $\star $ and $I,J\in F(D),~(IJ)^{\star }=(IJ^{\star })^{\star
}=(I^{\star }J^{\star })^{\star }$, (the $\star $-\emph{multiplication}),
(b) $I+^{\star }J=$ $(I+J)^{\star }$ $=(I+J^{\star })^{\star }=(I^{\star
}+J^{\star })^{\star }$ (the $\star $-\emph{sum}) and (c) $(I^{\star }\cap
J^{\star })^{\star }=I^{\star }\cap J^{\star }$($\star $\emph{-intersection}%
). A reader in need of a quick review of star operations may consult
sections 32 and 34 of \cite{Gil 1972}. For our purposes we include the
following.

To each star operation $\star $ we can associate a star operation $\star
_{s} $ defined by $I^{\star _{s}}=\bigcup \{\,J^{\star }\mid J\subseteq I$
and $J\in f(D)\,\}$. A star operation $\star $ is said to be of \emph{finite
type}, or of \emph{finite character}, if $I^{\star }=I^{\star _{s}}$ for all 
$I\in F(D)$. Indeed for each star operation $\star $, $\star _{s}$ is of
finite character. Thus if $\star $ is of finite character $I\in F(D)$ is a $%
\star $-ideal if and only if for each finitely generated subideal $J$ of $I$
we have $J^{\star }\subseteq I$. Also it is easy to see that for a nonzero
finitely generated ideal $I$ we have $I^{\star }=I^{\star _{s}}$. For $I\in
F(D)$, let $I_{d}=I$, $I^{-1}=D:_{K}I=\{\,x\in K\mid xI\subseteq D\,\}$, $%
I_{v}=(I^{-1})^{-1}$, $I_{t}=\bigcup \{\,J_{v}\mid J\subseteq I$ and $J\in
f(D)\,\}=I_{v_{s}}$ and so the $t$-\emph{operation} is an example of a star
operation of finite character. Star operations of finite character will
figure prominently in our discussions.

Let $I$ be anonzero integral ideal (i.e $I\subseteq D$) such that $I^{\star
}\neq D$, the definition of a $\star $-operation of finite type allows for
the existence of a \emph{maximal integral }$\star $\emph{-ideal} containing $%
I$, via Zorn's Lemma. A maximal $\star $-ideal can be shown to be a prime
ideal. (We say that $I$ is a \emph{proper }(or a\emph{\ proper }$\star $-%
\emph{ideal}) we shall mean that $I$ is an integral ideal properly contained
in $D$.)\emph{\ }Thus if $D$ is a domain that is not a field and if $\star $-%
$Max(D)$ denotes the set of maximal $\star $-ideals of $D$, for a star
operation $\star $ of finite character, then $\star $-$Max(D)\neq \phi $. A
fractional ideal $I$ is called $\star $\emph{-invertible} if $%
(II^{-1})^{\star }=D$. It is well known that if $I$ is $\star $-invertible
for a finite character star operation $\star $ then $I^{\star }$ and $I^{-1}$
are of finite type and that every $\star $-invertible $\star $-ideal is 
\emph{divisorial} (i.e. is a $v$-ideal).

An integral domain $D$ is a \emph{Prufer }$v$-\emph{multiplication domain
(PVMD)} if every nonzero finitely generated ideal of $D$ is $t$-invertible.
Denote the set of all $\star $-invertible fractional $\star $-ideals of $D$
by $Inv_{\star }(D)$ and note that $Inv_{\star }(D)$ is a group under $\star 
$-multiplication. We also note that $Inv_{\star }(D)$ is a partially ordered
group with order induced by $I\leq J$ $\Leftrightarrow $ $(JI^{-1})^{\star
}\subseteq D$. For a star operation $\star $ of finite character we show
that $Inv_{\star }(D)$ is a Riesz group if and only if every integral $\star 
$-invertible $\star $-ideal of $D$ is a primal element of $Inv_{\star }(D)$
under the induced order mentioned above. Given two star operations $\mu $
and $\rho $ we say that $\mu \leq \rho $, if $A^{\mu }\subseteq A^{\rho }$
for all $A\in F(D);$ equivalently $\mu \leq \rho $ if $(A^{\mu })^{\rho
}=(A^{\rho })^{\mu }=A^{\rho }$ for all $A\in F(D)$. Thus if $\mu \leq \rho $
every $\rho $-ideal is a $\mu $-ideal. Generally $\star \leq v$ for every
star operation $\star $.

Consider the set $I^{\star }(D)$ of integral $\star $-ideals of $D$ and note
that $I^{\star }(D)$ is a monoid under $\star $-multiplication: $(\times
^{\star })$, with $I\times ^{\star }J=($ $IJ)^{\star }$, with $D$ as its
identity. Also if we define the partial order $\leq $ on $F^{\star
}(D)=\{I^{\star }|I\in F(D)\}$ by $A\leq B$ if and only if $A\subseteq B$, $%
I^{\star }(D)$ is a lattice under $\star $-sum $+^{\star }$ defined by $%
I+^{\star }J=$ $(I+J)^{\star }$ and $\star $-intersection $\cap ^{\star }$
defined by $I\cap ^{\star }J=I\cap J$, with $\sup (I,J)=I+^{\star
}J=(I+J)^{\star }=I\vee J$ and $\inf (I,J)=I\cap ^{\star }J=$ $I\cap
J=I\wedge J$, for all $I,J\in F^{\star }(D)$. What is interesting is that
the order on $F^{\star }(D)$ is compatible with $\times ^{\star }$ and that
makes $F^{\star }(D)$ a lattice ordered monoid, with the property that $%
\times ^{\star }$ \emph{distributes} over $\vee _{\alpha }$ for all $\alpha $%
, in some index set $I$, i.e., $I\times ^{\star }(\vee _{\alpha \in
J}I_{\alpha })=\vee _{\alpha \in J}I\times ^{\star }I_{\alpha }$. It turns
out that the set $\mathcal{L}_{\star }(D)=I^{\star }(D)\cup \{0\}$ is a
multiplicative lattice with the least element $0$ and the greatest element $D
$, i.e., 
\begin{equation*}
\mathcal{L}_{\star }(D)=\langle I^{\star }(D)\cup \{0\},+^{\star },\times
^{\star },\vee ,\wedge ,D,0\rangle .
\end{equation*}%
This lattice is complete in that $\vee S=\vee \{s|s\in S\}$ and $\wedge
S=\wedge \{s|s\in S\}$ exist for each subset $S$ of $\mathcal{L}_{\star }(D)$%
. Because our work takes us to groups of divisibility, order on which is
defined by reverse containment, we are forced to consider defining $A\leq B$
if and only if $A\supseteq B$. With this definition of order we get 
\begin{equation*}
\Gamma _{\star }(D)=\langle I^{\star }(D)\cup \{0\},+^{\star },\times
^{\star },\wedge ,\vee ,0,D\rangle ,
\end{equation*}%
and with this definition of order, $0$ become the largest and $D$ the
smallest element of $\Gamma _{\star }(D)$. That's not all, $\sup (I,J)=I\cap
J$, $\inf (I,J)=I+^{\star }J$ and $\times ^{\star }$ distributes over $%
\wedge _{\alpha \in J}$, that is $I\times ^{\star }(\wedge _{\alpha \in
J}I_{\alpha })=\wedge _{\alpha \in J}I\times ^{\star }I_{\alpha }$ in $%
\Gamma _{\star }(D)$, if $\wedge _{\alpha \in J}I_{\alpha }\neq 0$. So,
though $\vee S$ and $\wedge S$ exist for each subset $S$ of $\Gamma _{\star
}(D)$, it's not quite a multiplicative lattice. Let's call $\Gamma _{\star
}(D)$ a \emph{reverse multiplicative lattice}. Because $0$ does not have a
function in matters of divisibility we shall consider $\Delta _{\star
}(D)=\Gamma _{\star }(D)\backslash \{0\}$, even at the cost of having to
settle for $\vee S$ existing for finite subsets $S$ of $I^{\star }(D)$. If
we restrict our attention to the monoid $\varphi _{\star }(D)$ of nonzero
integral $\star $-ideals of finite type under $\star $-multiplication with $%
I\leq J$ if and only if $I\supseteq J$, then $\varphi _{\star }(D)$ is a $%
\wedge $-semilattice with $I\times ^{\star }(J\wedge K)=I\times ^{\star
}J\wedge I\times ^{\star }K$.

Having described the set of star ideals of a domain, at such length, we owe
it to the reader to mention a recently studied generalization of Riesz
groups and hence of Riesz monoids. In \cite{YZ 2011} a directed p.o. group $%
G=\langle G,\ast ,e,\leq \rangle $ was called a pre-Riesz group if for $%
e<x_{1},x_{2},...,x_{n}\in G$, with $\mathcal{L(}\{x_{1},x_{2},...,x_{n}\})%
\neq \mathcal{L(}\{e\})$, there exists $r\in G$ such that $e<r\leq
x_{1},x_{2},...,x_{n}$. Here $\mathcal{L}(S)$ denotes the set of lower
bounds of a subset $S$ of $G$. Encouraged by our results on pre-Riesz groups
we defined and studied pre-Riesz monoids in \cite{LYZ 2020}. The definition
boils down, in the commutative case, to: a pre-Riesz monoid is a directed
p.o. monoid $M$, with least element $e$ such that for any $%
x_{1},x_{2},...,x_{n}\in M\backslash \{e\}$ $glb(x_{1},x_{2},...,x_{n})=e$
or there is $r\in M$ with $e<r\leq x_{1},x_{2},...,x_{n}$, where $glb$
stands for the greatest lower bound or $inf$.

Obviously if $G$ is a pre-Riesz group, then $G^{+}$ is a pre-Riesz monoid.
It turns out that some very interesting examples of pre-Riesz monoids are
afforded by $\Delta _{\star }(D)$ and its subsets. Being a lattice, with the
least element $D$ playing the role of the identity of the monoid we have:
for all 

\begin{center}
$D<A_{1},A_{2},...,A_{r}\in \Delta _{\star }(D)$, $\inf
(A_{1},A_{2},...,A_{r})\neq D$ $\Rightarrow $ $D<\inf
(A_{1},A_{2},...,A_{r})\leq A_{1},A_{2},...,A_{r}$.
\end{center}

Frankly, while a pre-Riesz monoid seemed to have plenty of interesting
features, $\Delta _{\star }(D)$ appears lackluster and devoid of teeth, as
it were. To give $\Delta _{\star }(D)$ \textquotedblleft
teeth\textquotedblright , designate a subset $\Pi $ of proper nonzero $\star 
$-ideals of $D$ and impose the condition that for all $%
D<A_{1},A_{2},...,A_{r}\in \Delta _{\star }(D)$, with $\inf
(A_{1},A_{2},...,A_{r})\neq D$ there is $A\in \Pi $ such that $A\leq
A_{1},A_{2},...,A_{r}$.

Now as $\inf (A_{1},A_{2},...,A_{r})$ means $\star $-sum of $A_{i}$ there
seems to be no reason to repeat the $A_{i}$. Thus we call the monoid $\Delta
_{\star }(D)$ a $\Pi $-Riesz monoid if for $A\in \Delta _{\star }(D),$ $%
A\neq D$ implies that there is $X\in \Pi $ such that $X\leq A$. We can
denote this setup as $D=(D,\Delta _{\star }(D),\Pi _{D})$ where $\Pi _{D}$
stands for $\Pi $ with reference to $D$. This designation gives rise to a
number of possibilities. For instance $\Pi _{D}$ can be the set of proper
nonzero principal (invertible, $t$-invertible $t$-ideals, $t$-ideals of
finite type, divisorial ideals) of $D$.

If we are no longer dealing with groups related to the group of divisibility
we can revert to the set up of $\mathcal{L}_{\star }(D)$ and our new setup
can be $D=(D,\mathcal{L}_{\star }(D)\backslash \{0\};\Pi _{D})$ defined by $%
A\neq D$ implies that there is $X\in \Pi $ such that $A\subseteq X$. Indeed
we can vary the star operation $\star $, setting $\star =d,t$ etc. and,
keeping the choices for $\Pi $ the same, we can vary $D$ as Noetherian or
some variation of it. We note that the $\wedge $-semilattice $\varphi
_{\star }(D)$ is a pre-Riesz monoid under $\star $-multiplication, if we set 
$\Pi _{D}=\varphi _{\star }(D)$. Also a p.o. monoid $M$ is called $\wedge $%
-smooth if whenever $x\wedge y$ exists in $M$ we have $z\ast (x\wedge
y)=z\ast x\wedge z\ast y$ and if $x\wedge (y\wedge z)$ or $(x\wedge y)\wedge
z$ exists then $x\wedge (y\wedge z)=(x\wedge y)\wedge z$.

One may ask whether Riesz monoids satisfy the Riesz interpolation, as do
Riesz groups. The answer is yes and can be readily verified as we shall show
below. Having shown that, in section \ref{S2}, we study the conditions under
which a Riesz monoid generates a Riesz group. Also in section \ref{S2}, we
study $\star $-Schreier domains as integral domains $D$ such that the group $%
Inv_{\star }(D)$ of $\star $-invertible $\star $-ideals of $D$, under $\star 
$-multiplication, is a Riesz group and show that $D$ is $\star $-Schreier if
and only if every integral $\star $-invertible $\star $-ideal of $D$ is a
primal element of $Inv_{\star }(D)$.

In section \ref{S3}, we study the $\Pi $-\emph{pre-Riesz monoids} in their
various guises and show for instance that a Noetherian domain is a Dedekind
domain (resp., a PID) if and only if $D=(D,\mathcal{L}_{d}(D),\Pi _{D}\}$
where $\Pi _{D}$ consists of proper invertible ideals (resp., nonzero
principal ideals) of $D$. Finally, in section \ref{S4}, we call an element $%
h $ in a $\wedge $-smooth pre-Riesz monoid $M$ a \emph{homogeneous element}
if $h>e$ and for all $e<u,v\leq h$ we have $e<t\leq u,v$. We show that in a $%
\wedge $-smooth pre-Riesz monoid $M$ a sum/product of finitely many
homogeneous elements can be expressed as a finite sum/product of mutually
disjoint homogeneous elements, up to some permutation of the constituents.

By calling a $\wedge $-smooth pre-Riesz monoid $M$ \emph{virtually factorial}
if $M$ is generated by its homogeneous elements and by observing that a
homogeneous element of the $\wedge $-semilattice $\varphi _{\star }(D)$ is
precisely a homogeneous ideal of $D$ we show: Let $\star $ be a finite
character star operation defined on an integral domain $D$. Then the monoid $%
\varphi _{\star }(D)$ of nonzero $\star $-ideals of finite type of $D$,
under $\star $-multiplication, is a $\wedge $-smooth pre-Riesz monoid with
order defined by $I\leq J$ if and only if $I\supseteq J$, for $I,J\in
\varphi _{\star }(D)$. Moreover the following hold: (1) $D$ is a $\star $-SH
domain if and only if $\varphi _{\star }(D)$ is a v-factorial monoid, (2) If 
$h(\varphi _{\star }(D))\neq \phi $, then $H_{h}(\varphi _{\star }(D))$ is a
v-factorial monoid. For $\star $-SH domains, the reader may consult section %
\ref{S4} or \cite{AZ 2019}.

\section{\textbf{Riesz monoids\label{S2}}}

We first show in this section that a Riesz monoid satisfies the $(2,2)$
Riesz interpolation and wave hands about the $(m,n)$ interpolation.

\begin{theorem}
\bigskip \label{Theorem QX} 
The following are equivalent for a commutative cancellation divisibility
monoid $M$. 

\begin{enumerate}
\item Every $x\in M^{+}$ is primal.

\item For all $a,b,x,y\in M^{+}$ with $a,b\leq $ $x,y$ there is $z$ such
that $a,b\leq z\leq x,y$.

\item For all $a,b,x_{1},$.$..,x_{n}\in M^{+}$ with $a,b\leq $ $x_{i}$ there
is $z$ such that $a,b\leq z\leq x_{i},$ $i=1,...,n$.

\item For all $a_{i},b_{j}\in M^{+}$ with $a_{i}\leq $ $b_{j}$ there is $z$
such that $a_{i}\leq z\leq b_{j}$.
\end{enumerate}
\end{theorem}

Part (2) of Theorem \ref{Theorem QX} is also called the $(2,2)$ Riesz
Interpolation Property. With some effort one can show that indeed a Riesz
monoid satisfies, and is characterized by, $(m,n)$ interpolation for all
positive integral $m$ and~$n$.

\begin{proof}
(1) $\Rightarrow $ (2). Suppose every positive element of $M$ is primal. Let 
$a,b\leq x,y.$ Then $x=x_{1}\ast a=x_{2}\ast b$ and $y=y_{1}\ast a=y_{2}\ast
b......................................(1)$

Since $x_{1}\ast a=x_{2}\ast b,$ $b\leq x_{1}\ast a$.

Also since $b$ is primal $b=b_{1}\ast b_{2}$ where $b_{1}\leq x_{1}$ and $%
b_{2}\leq a$........... (2)

Let $x_{1}=x_{1}^{\prime }\ast b_{1}$ and $a=a_{1}\ast b_{2}.$

Then $x_{1}\ast a=x_{2}\ast b$ can be written as $x_{1}^{\prime }\ast
b_{1}\ast a_{1}\ast b_{2}=x_{2}\ast b$, or $x_{1}^{\prime }\ast a_{1}\ast
b_{1}\ast b_{2}=x_{2}\ast b.$

Noting that $b=b_{1}\ast b_{2}$ and cancelling $b$ from both sides of the
previous equation we get

$x_{1}^{\prime }\ast a_{1}=x_{2}.$
.........................................................(3)

Since $a_{1}\ast b_{2}=a$ we have $a,b\leq a_{1}\ast
b.............................(4)$

Using the value of $x_{2}$ we have $a_{1}\ast b\leq x$ ...............(5)

(Note: $x=x_{2}\ast b=(x_{1}^{\prime }\ast a_{1})\ast b)$

Now consider $y_{1}\ast a=y_{2}\ast b.$

Using $a=a_{1}\ast b_{2}$ and $b=b_{1}\ast b_{2}$ we have $y_{1}\ast
a_{1}\ast b_{2}=y_{2}\ast b_{1}\ast b_{2}.$

Cancelling $b_{2}$ from both sides we get $y_{1}\ast a_{1}=y_{2}\ast b_{1}.$

So that $b_{1}\leq y_{1}\ast a_{1}$ and as $b_{1}$ is primal we have $%
b_{1}=b_{3}\ast b_{4}$ where $b_{3}\leq y_{1}$ and $b_{4}\leq a_{1}.$

Writing $y_{1}=y_{1}^{\prime }\ast b_{3}$ and $a_{1}=a_{1}^{\prime }\ast
b_{4}$ we can express $y_{1}\ast a_{1}=y_{2}\ast b_{1}$ as $y_{1}^{\prime
}\ast b_{3}\ast a_{1}^{\prime }\ast b_{4}=y_{2}\ast b_{1}.$ Cancelling $%
b_{1}=b_{3}\ast b_{4}$ from both sides we get $y_{2}=y_{1}^{\prime }\ast
a_{1}^{\prime }.$ This gives $y=y_{2}\ast b=y_{1}^{\prime }\ast
a_{1}^{\prime }\ast b=y_{1}\ast a.$ Now as $y_{1}^{\prime }\leq y_{1}$ we
get $y_{1}=y_{4}\ast y_{1}^{\prime }$ which on substituting in $%
y_{1}^{\prime }\ast a_{1}^{\prime }\ast b=y_{1}\ast a$ gives $y_{1}^{\prime
}\ast a_{1}^{\prime }\ast b=y_{4}\ast y_{1}^{\prime }\ast a$ and cancelling $%
y_{1}^{\prime }$ we get $y_{4}\ast a=a_{1}^{\prime }\ast b$ and so $a\leq
a_{1}^{\prime }\ast b.$ That is $a,b\leq a_{1}^{\prime }\ast b$ and $%
a_{1}^{\prime }\ast b\leq y.$ But as $a_{1}^{\prime }\leq a_{1}$ and $%
x_{2}=x_{1}^{\prime }\ast a_{1}$ we have $a_{1}^{\prime }\ast b\leq
x_{2}\ast b=x.$ So we have $z=a_{1}^{\prime }\ast b$ such that $a,b\leq
z\leq x,y.$

(2) $\Rightarrow $ (1). 
%
%
%
%
Let $a\leq b\ast c$. Then as $a,b\leq b\ast c,$ $a\ast b$ there is $x$ such that 
\[
a,b\leq x\leq b\ast c, a\ast b. \tag{i}
\]
Now as $a\leq x$ we have 
\[
x=x_{1}\ast a. \tag{ii} 
\]
Also as $b\leq x$ we have 
\[
x=x_{2}\ast b. \tag{iii}
\]
Using (i) and (iii), we have $x_{2}\leq a$ and $x_{2}\leq c.$ Now as $x_{2}\leq a,$
setting $a=x_{3}\ast x_{2}$ we have from $x_{1}\ast a=x_{2}\ast b$, the
equation $b=x_{1}\ast x_{3}.$ So $a\leq b\ast c$ implies that $a=x_{2}\ast
x_{3},$ with $x_{2},x_{3}\in M^{+}$ such that $x_{3}\leq b$ and $x_{2}\leq
c. $


((2) $\Rightarrow $ (3). Let $a,b\leq $ $x_{i},$ for $i=1,...,n,$ in $M.$ If 
$n=1$ we have $z=x_{1}$ and if $n=2$ then by (2) there is a $z\in M$ such
that $a,b\leq z\leq $ $x_{1},x_{2}.$ Suppose that for all $2\leq r\leq n$ we
have a $z_{r}$ such that $a,b\leq z_{r}\leq $ $x_{1},x_{2},...,x_{r}$ and
consider $a,b\leq $ $x_{1},x_{2},...,x_{n+1}.$ By the induction hypothesis
there is $z_{n}$ such that $a,b\leq z_{n}\leq $ $x_{1},x_{2},...,x_{n}$ and
this gives $a,b\leq z_{n},x_{n+1}.$ Now by (2) we have a $z_{n+1}$ such that 
$a,b\leq z_{n+1}\leq z_{n},x_{n+1}.$ But since $z_{n+1}\leq z_{n}\leq $ $%
x_{1},x_{2},...,x_{n}$ we have the result that $a,b\leq z_{n+1}\leq $ $%
x_{1},x_{2},...,x_{n+1}.$ Thus, by induction, if $a,b\leq $ $x_{i}$ for $%
i=1,...,n,$ in $M$, then there is $z$ in $M$ such that $a,b\leq z\leq $ $%
x_{i}$ for $i=1,...,n.$

(3) $\Rightarrow $ (4). Suppose $a_{i},b_{j}\in M^{+}$ with $a_{i}\leq $ $%
b_{j}$ for $i=1,...,m,$ $j=1,...,n.$ Let $m=1.$ Then $z=a_{1}$ will do the
job and if $m=2,$ the statement of (3) gives a $z$ such that $%
a_{1},a_{2}\leq $ $z\leq b_{1},b_{2},...,b_{n}.$ Now suppose that for all $%
2\leq r\leq m$ we have that $a_{1},...,a_{r}\leq b_{1},b_{2},...,b_{n}$
implies the existence of $z_{r}$ with $a_{1},...,a_{r}\leq z_{r}\leq
b_{1},b_{2},...,b_{n}.$ Consider $a_{1},...,a_{r},a_{r+1}\leq
b_{1},b_{2},...,b_{n}$. But as there is $z_{r},$ by induction hypothesis, we
have $z_{r},a_{r+1}\leq b_{1},b_{2},...,b_{n}$ and (3) applies to give $%
z_{r+1}$ such that $z_{r},a_{r+1}\leq z_{r+1}\leq b_{1},b_{2},...,b_{n},$
which completes the job because $a_{1},...,a_{r}\leq z_{r}.$

(4) $\Rightarrow $ (2). Obvious once we put $m=n=2.)$
\end{proof}

Call a subset $S$ of a monoid $M$ \emph{conic} if $x\ast y=e$ implies $x=e=y$%
, for all $x,y\in S$. In a p.o. group $G$ the sets $G^{+}$ and $-G^{+}$ are
conic. If $D$ is an integral domain then the set $m(D)$ of nonzero principal
ideals of $D$ is a monoid under multiplication, with identity $D$, ordered
by $aD\leq bD$ $\Leftrightarrow $ there is $c\in D$ such that $%
bD=acD\Leftrightarrow aD\supseteq bD$. The monoid $m(D)$ is cancellative too
and in $m(D)$ $xDyD=D\Rightarrow xD=yD=D$. So, $m(D)$ is a divisibility
cancellative conic monoid. The monoid $m(D)$ is of interest because of the
manner in which it generates a group. We know how the field of quotients of
a domain is formed as a set of ordered pairs, each pair representing an
equivalence class with $(a,b)=(c,d)$ $\Leftrightarrow da=bc$ and then we
represent the pair $(a,b),$ $b\in D\backslash \{0\}$ by $\frac{a}{b}=ab^{-1}$%
.

Now the group of $m(D)$ gets the form

\begin{center}
$G(D)=\{\frac{a}{b}D|\frac{a}{b}\in qf(D)\backslash \{0\}\}$, ordered by $%
\frac{a}{b}D\leq \frac{c}{d}D$ $\Leftrightarrow \frac{a}{b}D\supseteq \frac{c%
}{d}D$

$\Leftrightarrow $

there is $hD\in m(D)$ such that $\frac{a}{b}DhD=\frac{a}{b}D$, so that $m(D)$
is the positive cone of $G(D)$.

The group $G(D)$ gets the name \emph{group of divisibility} of $D$ (actually
of $m(D))$. Now any divisibility monoid that is also a cancellative and
conic monoid $M$, with least element $e$ can be put through a similar
process of forming equivalence classes of ordered pairs to get its group of
divisibility like group $G(M)=\{a\ast b^{-1}|a,b\in M\}$ with $x\leq y$ in $%
G(M)$ $\Leftrightarrow x\ast h=y$ for some $h\in M$. (Here $a\ast b^{-1}$
becomes $a-b$, if $\ast =+$.$)$
\end{center}

\begin{corollary}
\label{Corollary RX}A Riesz Monoid $M$ has the pre-Riesz property. Also $%
M^{+}$ is conic for a Riesz monoid $M$.
\end{corollary}

\begin{proof}
Let $e\leq x,y$ in $M$ and suppose that there is $g\in M$ such that $g$ is
not greater than or equal to $e$ yet $g\leq x,y,$ that is $e,g\leq x,y.$
Then by the $(2,2)$ interpolation property there is $r\in M$ such that $%
e,g\leq r\leq x,y.$ But then $r>e,$ as $r\geq e$ and $r\neq e$ because $%
r\geq g.$ Next suppose $x,y\geq e$. If $x\ast y=e$ and say $x\neq e,$ then
we have $e,x\leq x,x\ast y$ and by the $(2,2)$ interpolation there is $r$
such that $e<r\leq x,x\ast y$ contradicting the fact that $x\ast y=e.$
\end{proof}

Well a p.o. monoid $M$ is a p.o. group if every element of $M$ has an
inverse and obviously if a p.o. monoid is a Riesz monoid and a group, it is
a Riesz group. This brings up the question: Let $M$ be a Riesz monoid and $%
M^{+}$ the positive cone of it, will $M^{+}$ generate a Riesz group? As we
shall be mostly concerned with monoids $M$ with $e$ the least element, i.e. $%
M=M^{+}$, we remodel the question as: Let $M$ be a Riesz monoid with $%
M=M^{+} $ the positive cone of it, will $M$ generate a Riesz group? The
following result whose proof was indicated to me by G.M. Bergman, in an
email, provides the answer. (To preserve an example of his frank
conversational style, I haven't changed the arguments much. By the way,
George has always been a source of inspiration for me. Let me put it this
way, I wrote my doctoral thesis, looking at and at times looking into the
big book called \textquotedblleft Bergman's thesis" on the table in the
graduate study room at Bedford College London.)

\begin{theorem}
\label{Theorem SX}Suppose $M$ is a cancellative abelian monoid, which is
\textquotedblleft conical\textquotedblright , i.e., no two non identity
elements add/multiply to $e$, and which we partially order by divisibility;
and suppose every element of $M$ is primal, namely, that with respect to the
divisibility induced order, 
\begin{equation*}
x\leq a\ast b\implies x=u\ast v,\text{ for some $u\leq a$ and $v\leq b$
....(1)}
\end{equation*}%
%
%
%
%
%
%
%
%
%
%
%
%
%
%
%
%
%
%
%
%
Then the group generated by $M$ is a Riesz group.
\end{theorem}

\begin{proof}
Let us rewrite (1) by translating all the inequalities into their
divisibility statements; so that $x\leq a\ast b$ becomes $x\ast y=a\ast b$
for some $y$ and $u\leq a$ becomes $a=u\ast u_{1}$, and similarly for the
last inequality; and finally, let us rename the elements more
systematically; in particular, using $a,b,c,d$ for the above $x,y,a,b$. Then
we find that (1) becomes $a\ast b=c\ast d\Rightarrow a=a_{1}\ast
a_{2},c=a_{1}\ast b_{1},d=a_{2}\ast b_{2}$ for some $a_{1},a_{2},b_{1},b_{2}%
\in M.$ Now if we substitute the three equations to the right of the "$%
\Rightarrow $" into the equation before the "$\Rightarrow $", and use
cancellativity, we find that $b=b_{1}\ast b_{2}$; so the full statement is
(2) $a\ast b=c\ast d\Rightarrow a=a_{1}\ast a_{2},~b=b_{1}\ast
b_{2},~c=a_{1}\ast b_{1},~d=a_{2}\ast b_{2},$ for some $%
a_{1},a_{2},b_{1},b_{2}\in M$. Now let $G$ be the group generated by $M$,
ordered so that $M$ is the positive cone. We want to show $G$ has the Riesz
Interpolation Property. So suppose that in $G$ we have $p,q\leq r,s$. We can
write these inequalities as (3) $r=p\ast a,s=p\ast c,r=q\ast d,s=q\ast b$
where $a,b,c,d\in M$. Now the sum of the first and last equations gives a
formula for $r\ast s$, and so does the sum of the second and third
equations. Equating the results, and cancelling the summands $p\ast q$ on
each side, we get an equation in $M:a\ast b=c\ast d$. Hence we can apply (2)
to get decompositions of $a,b,c,d$, and substitute these into (3), getting
(4) $r=p\ast a_{1}\ast a_{2},s=p\ast a_{1}\ast b_{1},r=q\ast a_{2}\ast
b_{2},s=q\ast b_{1}\ast b_{2}$. Equating the first and third equations (or
if we prefer, the second and fourth) and cancelling the common term $a_{2}$
(respectively, the common term $b_{1}$), we get (whichever choice we have
made) (5) $p\ast a_{1}=q\ast b_{2}=t$. The element given by (5) is clearly $%
\geq p,q$, while from (4) (using whichever of the equations for $r$ we
prefer and whichever of the equations for $s$ we prefer, say $r=p\ast
a_{1}\ast a_{2}$ (a) and $s=q\ast b_{1}\ast b_{2}$ (b), using $q\ast b_{2}=$ 
$p\ast a_{1}=t$ ), we see that it is $\leq r,s$. So this is the element
whose existence is required for the ($(2,2)$) Riesz interpolation property
for $G$.
\end{proof}

Let $\mathcal{I}_{\star }(D)$ be the set of integral $\star $-invertible $%
\star $-ideals and note that $\mathcal{I}_{\star }(D)$ is a monoid under $%
\star $-multiplication. Note also that $\mathcal{I}_{\star }(D)$ is
partially ordered by $I\leq J$ if and only if $I\supseteq J$. Indeed $%
J\subseteq I$ if and only if $(JI^{-1})^{\star }=H\subseteq D$, if and only
if $J=(IH)^{\star }$, and as $J,$ $I$ are $\star $-invertible, $H$ is $\star 
$-invertible and integral. Thus in $\mathcal{I}_{\star }(D)$, $I\leq J$ $%
\Leftrightarrow J=(IH)^{\star }$ for some $H\in \mathcal{I}_{\star }(D)$. In
other words $\mathcal{I}_{\star }(D)$ is a divisibility p.o. monoid. Because 
$\mathcal{I}_{\star }(D)$ involves only $\star $-invertible $\star $-ideals,
it is cancellative too. Finally $\mathcal{I}_{\star }(D)$ is directed
because of the definition of order. That $Inv_{\star }(D)$ is generated by $%
\mathcal{I}_{\star }(D)$ follows from the fact that every fractionary ideal
of $D$ can be written in the form $A/d$ where $A\in F(D)$ and $d\in
D\backslash \{0\}$. Finally, the partial order in $Inv_{\star }(D)$ gets
induced by $\mathcal{I}_{\star }(D)$ in that for $I,J\in $ $Inv_{\star }(D)$
we have

\begin{center}
$I\leq J$ $\Leftrightarrow J\subseteq I\Leftrightarrow (JI^{-1})^{\star }\in 
$ $\mathcal{I}_{\star }(D)$.
\end{center}

Call $I\in \mathcal{I}_{\star }(D)$ $\star $-primal if for all $J,K\in 
\mathcal{I}_{\star }(D)$ $I\leq (JK)^{\star }$ we have $I=(I_{1}I_{2})^{%
\star }$ where $I_{1}^{\star }\leq J$ and $I_{2}^{\star }\leq K$. Call $D$ $%
\star $-Schreier, for a star operation $\star $ of finite character, if
every integral $\star $-invertible $\star $-ideal of $D$ is primal.

\begin{remark}
\label{Remark SX1}For the group-theoretic version of Theorem \ref{Theorem SX}%
, see (4) of Theorem 2.2 of \cite{Fu 1965}. However, Theorem \ref{Theorem SX}
is not a repeat in that the environment is different. For instance, here we
had to prove that a Riesz monoid is conic, while in the p.o. group
environment $G^{+}$ being conic is a given. Let's also note here that $%
M=\{x\in Q^{+}|x=0$ or $x\geq 1\}$ is a monoid under addition with $%
G(M)=\{x-y|x,y\in M\}=Q$ a patently Riesz group, being totally ordered. But
as we shall check below $M$ is not a Riesz monoid, because $M$ is not a
divisibility monoid under the order induced on it by $Q$ and so $M$ cannot
serve as a positive cone. The reason is, say, $1\leq 1$.$1$ in $M$ under the
order induced by $G(M)=Q$, but $.1\notin M$. So the natural order on $M$
iduced by divisibility is different from the order induced by $Q$ and $M$
cannot serve as a positive cone of the Riesz group $Q$. In other words when
a group $G$ generated by a monoid $M$ turns out to be Riesz we cannot
conclude that $M$ is a Riesz monoid without checking that $M=G^{+}$, under
the order induced by divisibility.
\end{remark}

\begin{proposition}
\label{Proposition TX}Let $\star $ be a finite character star operation
defined on $D$. Then $D$ is a $\star $-Schreier domain if and only if $%
Inv_{\star }(D)$ is a Riesz group under $\star $-multiplication and order
defined by $A\leq B\Leftrightarrow A\supseteq B$.
\end{proposition}

\begin{proof}
Suppose that $D$ is $\star $-Schreier, as defined above. That is each $I\in $
$\mathcal{I}_{\star }(D)$ is primal. The notion of $\star $-Schreier
suggests that we define $\leq $ by $A\leq B$ $\Leftrightarrow A\supseteq B.$
Then as for each pair of integral ideals $I,J$, $(IJ)^{\star }=D$ $%
\Rightarrow J^{\star }=I^{\star }=D,$ the same holds for members of $%
\mathcal{I}_{\star }(D)$ which are all integral $\star $-ideals. So $%
(IJ)^{\star }=D\Rightarrow I=J=D.$ and so $\mathcal{I}_{\star }(D)$ is
conic. Of course $\mathcal{I}_{\star }(D)$ is cancellative by the choice of
ideals and by the definition of order $\mathcal{I}_{\star }(D)$ is a
divisibility monoid. So by Theorem \ref{Theorem SX} $\mathcal{I}_{\star }(D)$
generates a Riesz group and by the above considerations $Inv_{\star }(D)$ is
generated by $\mathcal{I}_{\star }(D).$ Consequently $Inv_{\star }(D)$ is a
Riesz group. Conversely if $Inv_{\star }(D)$ is a Riesz group, with that
order defined on it, then $\mathcal{I}_{\star }(D)$ is the positive cone of
the Riesz group $Inv_{\star }(D)$ and so each element of $\mathcal{I}_{\star
}(D)$ must be primal.
\end{proof}

Proposition \ref{Proposition TX} brings together a number of notions studied
at different times. The first was \emph{quasi-Schreier}, study started in 
\cite{DM 2003} and completed in \cite{ADZ 2007}. The target in these papers
was studying $\mathcal{I}_{d}(D)$, i.e. the monoid of invertible integral
ideals of $D$, when $Inv_{d}(D)$ is a Riesz group. Another study targeting $%
\mathcal{I}_{t}(D)$, i.e. the monoid of $t$-invertible integral $t$-ideals
of $D$, for study along the same lines as above appeared in \cite{DZ 2011},
under the name $t$-\emph{Schreier domains.}

Now let's step back and require that every $\star $-invertible $\star $%
-ideal of $D$ be principal. Then in Proposition \ref{Proposition TX}, $%
\mathcal{I}_{\star }(D)$ is the monoid of nonzero principal integral ideals,
each of which is primal and the Riesz group $Inv_{\star }(D)$ consists just
of principal fractional ideals of $D$, and hence the group of divisibility
of $D$. But then $D$ is what was dubbed as a pre-Schreier domain in \cite%
{Zaf 1987}. Let us finally note that an integrally closed pre-Schreier
domain was initially introduced by Cohn in \cite{Coh 1968}, as a Schreier
ring. Cohn did say in \cite{Coh 1968} that the group of divisibility of a
Schreier domain was a Riesz group.

Now, looking at it from this angle, one can say that $D$ is a Schreier
domain if $D\backslash \{0\}$ is a multiplicative Riesz monoid. (Actually,
in view of our remarks abve, it would be better to say that $D$ is
pre-Schreier if $\{xD|x\in D\backslash \{0\}$ is Riesz monoid.) But the
story does not end here. In \cite{Coh 1968}, Cohn proved that if $D$ is a
Schreier ring and $X$ an indeterminate over $D$, then so is $D[X]$. Later,
there was a lot of activity on monoid domains with Matsuda, at the
forefront, converting a lot of results on polynomial rings to monoid rings.
He, more than, translated Cohn's result into the language of monoid rings by
saying in \cite{Mat 1988} that $D[X,S]$ is a Schreier domain if and only if $%
D$ and $K[X;S]$ are Schreier domains where $S$ is an integrally closed
additive Schreier semigroup (actually monoid), the idea of an additive
(pre-) Schreier monoid was born.

Let me mention for those who need to refresh their memory that a monoid $M$
is integrally closed if for $h\in q(M)=\{m_{1}\ast m_{2}^{-1}|m_{i}\in M\}$,
and for $n$ a natural number, $nh\in M^{+}$ implies $h\in M^{+}$, here $nh$
stands, for want of a better notation, for $h\ast h...n$ times. In other
words, $M$ is integrally closed if, whenever $nx\geq ny$ for some positive
integer $n$ we have $x\geq y$. Now it so happens that the description of an
additive (pre-) Schreier monoid provided in \cite{Mat 1988} matches our
description of Riesz monoids. Hence the general treatment of Riesz monoids,
using $\ast $ for the binary operation, so that additive Riesz monoids are
well represented.

Finally a word about monoid domains. Given a commutative additive monoid $M$
and a commutative ring $R$ we can form a set $R[X;M]=\{%
\sum_{i=1}^{k}r_{i}X^{m_{i}}|r_{i}\in R$ and $m_{i}\in M\}$. Defining
addition and multiplication as we do for polynomial addition and
multiplication using $X^{m_{1}}X^{m_{2}}=X^{m_{1}+m_{2}}$ we can make $%
R[X;M] $ a ring. Now according to Theorem 8.1 of Gilmer's book \cite{Gil
1984} $D[X;M]$ is a domain if and only if $D$ is a domain and $M$ is a
grading (additive, cancellative and torsionfree) monoid.

\section{Pre-Riesz monoids\label{S3}}

\begin{center}
If we adopt the language of $\mathcal{L}_{\star }(D)$ that takes $A\leq B$
if and only if $A\subseteq B$, the p.o. monoid $\mathcal{L}_{\star
}(D)\backslash \{0\}$ is a pre-Riesz monoid because $A_{1},A_{2},...,A_{r}%
\subseteq \sup (A_{1},A_{2},...,A_{r})\subsetneq D$, for all

$A_{1},A_{2},...,A_{r}\in \mathcal{L}_{\star }(D)\backslash \{0\}$ with $%
\sup (A_{1},A_{2},...,A_{r})\neq D$.

Further if we relax it to: for all 
\begin{eqnarray*}
A_{1},..A_{r} &\in &\mathcal{L}_{\star }(D)\backslash
\{0\},(A_{1},A_{2},...,A_{r})^{\star }\neq D \\
&\Rightarrow & \\
\exists D &\neq &A\in \mathcal{L}_{\star }(D)\backslash \{0\}
\end{eqnarray*}
such that $A_{i}\subseteq A$ we still get a pre-Riesz monoid. If this
doesn't sound very exciting let us designate a non-empty subset $\Pi _{D}$
of $\mathcal{L}_{\star }(D)$, say $\Pi _{D}$ is the set of proper nonzero
principal ideals of $D$. Now impose on $\mathcal{L}_{\star }(D)\backslash
\{0\}$ the condition: 
\begin{equation*}
\forall A\in \mathcal{L}_{\star }(D)\backslash \{0\}A\neq D\Rightarrow
\exists \pi \in \Pi _{D}
\end{equation*}
with $A\subseteq \pi $. Denote this domain by $D=(D,\mathcal{L}_{\star
}(D),\Pi _{D})$. Clearly this satisfies: 
\begin{equation*}
\forall A_{1},..A_{r}\in \mathcal{L}_{\star }(D)\backslash
\{0\},(A_{1},A_{2},...,A_{r})^{\star }\neq D\Rightarrow \exists \pi \in \Pi
_{D}
\end{equation*}
such that $A_{i}\subseteq \pi $
\end{center}

and so $\mathcal{L}_{\star }(D)\backslash \{0\}$ is a pre-Riesz monoid.

Now, as we have seen in the introduction above, this translates to: every
proper $\star $-ideal of $D$ is contained in a proper principal ideal of $D$%
. In other words, if $A$ is a $\star $-ideal that is not contained in any
principal ideal then $A=D$. This too may not be very exciting at a first
glance, even though every maximal $\star $-ideal of $D$ is principal as a
result. Next, take a pair of coprime elements $x,y$ of $D$. Then $%
(x,y)^{\star }=D$. Thus, at least, for every pair $x,y$ of coprime elements
in $D$ we have $(x,y)^{\star }=D$. Now let $d$ be an irreducible element of $%
D$ and suppose that $d|ab$ for some $a,b\in D$. If $d\nmid a$ and $d\nmid b$%
, then $D=((d,a)^{\star }(d,b)^{\star })^{\star }=(d^{2},da,db,ab)^{\star
}\subseteq dD$ a contradiction, because $d|ab$.

Thus an irreducible element is a prime in the domain 
\begin{equation*}
D=(D,\mathcal{L}_{\star }(D)\backslash \{0\},\Pi )
\end{equation*}%
,

for any star operation $\star $. Now for $\star =d$, the identity operation.
Next 
\begin{equation*}
D=(D,\mathcal{L}_{d}(D)\backslash \{0\},\Pi )
\end{equation*}

is a domain in which every proper nonzero ideal is contained in a principal
ideal, something stronger than what Cohn \cite{Coh 1968} called a pre-Bezout
domain. This domain 
\begin{equation*}
D=(D,\mathcal{L}_{d}(D)\backslash \{0\},\Pi )
\end{equation*}%
is even stronger than what was called a special pre-Bezout, or spre-Bezout
domain in \cite{DZ 2010}. Similarly if $D=(D,\mathcal{L}_{v}(D),\Pi )$, then 
$D$ is something stronger than a PSP-domain (every primitive polynomial over 
$D$ is super-primitive), also discussed in \cite{DZ 2010}. Can we find
domains that satisfy these properties? Yes indeed!

\begin{example}
\label{Example TX0} Let $Z,Q$ denote the ring of integers and its quotient
field respectively and let $X$ be an indeterminate over $Q$, then the ring $%
D=Z+XQ[X]$ is $(D,\mathcal{L}_{d}(D)\backslash \{0\},\Pi )$ where $\Pi $ is
the set of proper principal nonzero ideals of the domain $D$. 
\end{example}


\begin{illustration}
According to \cite[Theorem 4.21]{CMZ 1978} the prime ideals of $D$ are of
the form $pZ+XQ[X]$,$XQ[X]$ and height one principal primes of the form $%
f(X)D$ where $f(X)$ is irreducible in $K[X]$ and $f(0)=1$. Also, according
to \cite[Proposition 4.12]{CMZ 1978}, a general ideal of $D$ is of the form $%
I=f(X)(F+XQ[X])$ where $F$ is a $Z$-submodule of $Q$ such that $%
f(0)F\subseteq Z$. If $f(0)=0$, $f(X)$ is of the form $Xg(X)$ where $g(X)\in
Q[X]$ and so for each $z\in Z$ we have $I\subseteq zZ.$ If on the other hand 
$f(0)=d\neq 0$, then as $f(0)F\subseteq Z$ we have that $dF$ is an ideal of $%
Z.$ So $I=f_{1}(X)d(F+XQ[X])$ where $f_{1}(X)\in D$ such that $f_{1}(0)=1$
and $I=f_{1}(X)(dF+XQ[X]).$ If $dF$ is a proper ideal we need look no
further for proper principal ideals containing $I$. But if $dF=Z$ and $I$ is
proper then $f_{1}(X)$ must be variable and a product of powers of primes of 
$D$. This \textquotedblleft analysis\textquotedblright\ also establishes
that if no principal ideals of $D$ contain $I=f(X)(F+XQ[X])$, then $I$ must
be equal to $D$.
\end{illustration}

Now enough with this motivational chit chat. Let's state/prove some results.

A number of statements can be made in this connection, some examples are
given below. Yet considering the situation that we are faced up with, some
observations are in order.


\begin{observation}
\label{Observation TX1}Suppose that $P=\langle P,\leq \rangle $ is a non
empty poset with the properties that every element of $P$ precedes a maximal
element of $P$ and suppose that a non empty set $\Pi $ is chosen from $P$ by
some rule. Then every maximal element of $P$ is in $\Pi $ if and only if
every element of $P$ is required to precede some element of $\Pi $.

This, somewhat simple observation may, in some instances, have some
interesting consequences.
\end{observation}

\begin{lemma}
\begin{enumerate}
\item \label{Lemma TX2} If $\star \leq \rho $ are two star operations, where 
$\star $ is of finite type, and if every maximal $\star $-ideal of $D$ is a $%
\rho $-ideal, then the set $\Pi _{D}(\rho $-$inv))$ of $\rho $-invertible $%
\rho $-ideals coincides with the set $\Pi _{D}(\star $-$inv))$ of $\star $%
-invertible $\star $-ideals of $D$.

\item Suppose that $\star \leq \rho $ are two star operations, where $\star $
is of finite type. Then 
\begin{equation*}
D=(D,\mathcal{L}_{\star }(D)\backslash \{0\},\Pi _{D}(\rho -inv))
\end{equation*}
for $\Pi _{D}(\rho $-$inv)$ the set of proper $\rho $-invertible $\rho $%
-ideals of $D$, if and only if every maximal $\star $-ideal of $D$ is a $%
\rho $-invertible $\rho $-ideal of $D$.
\end{enumerate}
\end{lemma}

\begin{proof}
(1) Suppose every maximal $\star $-ideal is a $\rho $-ideal and let $I$ be a 
$\rho $-invertible ideal. Claim that $I$ is $\star $-invertible. For if not,
then $(II^{-1})^{\star }\subseteq P$ for some maximal $\star $-ideal $P$.
But then, as $P$ is a $\rho $-ideal, $((II^{-1})^{\star })^{\rho }\subseteq
P.$ Also as $((II^{-1})^{\star })^{\rho }=(II^{-1})^{\rho }=D$ because $\rho
\geq \star ,$ a contradiction. Next as a $\rho $-ideal is a $\star $-ideal
we conclude that a $\rho $-invertible $\rho $-ideal is a $\star $-invertible 
$\star $-ideal. Thus $\Pi _{D}(\rho $-$inv))\subseteq \Pi _{D}(\star $-$%
inv)).$ For the reverse containment note that as $\star \leq \rho ,$ $%
(II^{-1})^{\star }=D\Rightarrow (II^{-1})^{\rho }=D$, i.e, $\star $%
-invertible is $\rho $-invertible. Also as a $\star $-invertible $\star $%
-ideal is a $v$-ideal and hence a $\rho $-ideal, we have the conclusion.

(2) Suppose $D=(D,\mathcal{L}_{\star }(D)\backslash \{0\},\Pi _{D}(\rho $-$%
inv)).$ Then every maximal $\star $-ideal $M$ of $D$ is contained in a $\rho 
$-invertible $\rho $-ideal of $D$ and hence must be a $\rho $-invertible $%
\rho $-ideal. Conversely suppose that every maximal $\star $-ideal of $D$ is
a $\rho $-invertible $\rho $-ideal$.$ By (1), $\Pi _{D}(\rho $-$inv)=\Pi
_{D}(\star $-$inv)$. Thus $D=(D,\mathcal{L}_{\star }(D)\backslash \{0\},\Pi
_{D}(\rho $-$inv)).$
\end{proof}

\begin{proposition}
\label{Proposition UX0}

\begin{enumerate}
\item Let $\star >d$ be a finite type star operation defined on $D$. Then $%
D=(D,\mathcal{L}_{d}(D)\backslash \{0\},\Pi _{D})$ for $\Pi _{D}$ the set of
proper nonzero principal ideals (resp., $\star $-invertible $\star $-ideals, 
$\star $-ideals of finite type, $\star $-ideals, divisorial ideals) if and
only if every maximal ideal of $D$ is a principal ideal (resp., invertible
ideal, $\star $-ideal of finite type, divisorial ideal) of $D$.

\item An integral domain $D$ is $D=(D,\mathcal{L}_{d}(D)\backslash \{0\},\Pi
_{D})$ for $\Pi _{D}$ the set of proper nonzero principal ideals (resp., $t$%
-invertible $t$-ideals, $t$-ideals, $t$-ideals of finite type, divisorial
ideals) if and only if every maximal ideal of $D$ is a principal ideal
(resp., invertible ideal,, $t$-ideal, $t$-ideal of finite type, divisorial
ideal),

\item An integral domain $D$ is $D=(D,\mathcal{L}_{t}(D),\Pi _{D})$ for $\Pi
_{D}$ the set of proper nonzero principal ideals (resp.,invertible ideals, $%
t $-invertible $t$-ideals, $t$-ideals of finite type, divisorial ideals) if
and only if every maximal $t$-ideal of $D$ is a principal ideal (resp.,
invertible ideal, $t$-invertible $t$-ideal, $t$-ideal of finite type,
divisorial ideal) of $D$.
\end{enumerate}
\end{proposition}

\begin{proof}
In the presence of Observation \ref{Observation TX1} and Lemma \ref{Lemma
TX2}, it appears totally unnecessary to repeat the arguments required for
the proofs of (1) and (2). However, we take up selective cases from (1).
Suppose $D=(D,\mathcal{L}_{d}(D)\backslash \{0\},\Pi _{D})$ for $\Pi _{D}$
the set of proper $\star $-ideals ($\star $-ideals of finite type,
divisorial ideals). Then every maximal ideal is a $\star $-ideal ($\star $%
-ideal of finite type, divisorial ideal), by the condition. So $%
Max(D)\subseteq \Pi _{D}$ the set of proper $\star $-ideals ($\star $-ideals
of finite type, divisorial ideals). Whence we can say that every proper
ideal of $D$ is contained in a proper $\star $-ideal ($\star $-ideal of
finite type, divisorial ideal) of $D.$ Now (2) is just a special case of
(1). For (3), let $D=(D,\mathcal{L}_{t}(D),\Pi _{D})$ for $\Pi _{D}$ the set
of proper nonzero principal ideals (resp.,invertible ideals, $t$-invertible $%
t$-ideals, $t$-ideals of finite type, divisorial ideals). Then, by the
condition every maximal $t$-ideal of $D$ is principal (resp., invertible, a $%
t$-invertible $t$-ideal, a $t$-ideal of finite type, divisorial ideal). Thus
maximal $t$-ideals of $D$ are contained in $\Pi _{D}$ the set of proper
nonzero principal ideals (resp.,invertible ideals, $t$-invertible $t$%
-ideals, $t$-ideals of finite type, divisorial ideals) and one can say that
every proper $t$-ideal is contained in a nonzero proper principal ideal
(resp.,invertible ideal, $t$-invertible $t$-ideal, $t$-ideal of finite type,
divisorial ideal)
\end{proof}

Note that in case of (2) every maximal ideal being a $t$-ideal of finite
type ensures that every maximal $t$-ideal of $D$ is actually a maximal
ideal. Indeed if we suppose that $P$ is a maximal $t$-ideal that is not
maximal then $P$ is contained in a maximal ideal, say $M$, but $M$ is
already a $t$-ideal.

When it comes to the usual extensions of domains

\begin{center}
$D=(D,\mathcal{L}_{d}(D)\backslash \{0\},\Pi _{D})$ for suitable $\Pi _{D}$s
\end{center}

present interesting scenarios. We restrict to the star operations that are
defined for extensions considered.

\begin{proposition}
\label{Proposition UXA}

\begin{enumerate}
\item Let $D=(D,\mathcal{L}_{d}(D)\backslash \{0\},\Pi _{D})$ for $\Pi _{D}$
the set of proper nonzero principal ideals (resp., $t$-invertible $t$%
-ideals, $t$-ideals, $t$-ideals of finite type, divisorial ideals), let $X$
be an indeterminate over $D$ and let $R=D[X]$. Then it never is the case
that $R=(R,\mathcal{L}_{d}(R)\backslash \{0\},\Pi _{R})$ for $\Pi _{R}$ the
set of proper nonzero principal ideals (resp., $t$-invertible $t$-ideals, $t$%
-ideals, $t$-ideals of finite type, divisorial ideals) and

\item Let $D=(D,\mathcal{L}_{t}(D),\Pi _{D})$ for $\Pi _{D}$ the set of
proper nonzero principal ideals (resp., $t$-invertible $t$-ideals, $t$%
-ideals of finite type, divisorial ideals), let $X$ be an indeterminate over 
$D$ and let $R=D[X]$. Then $R=(R,\mathcal{L}_{t}(R)\backslash \{0\},\Pi
_{R}) $ for $\Pi _{R}$ the set of proper nonzero principal ideals (resp., $t$%
-invertible $t$-ideals, $t$-ideals of finite type, divisorial ideals) and
conversely.
\end{enumerate}
\end{proposition}

\begin{proof}
(1) Let $D=(D,\mathcal{L}_{d}(D)\backslash \{0\},\Pi _{D})$ for $\Pi _{D}$
the set of proper nonzero principal ideals (resp., $t$-invertible $t$%
-ideals, $t$-ideals, $t$-ideals of finite type, divisorial ideals). Then
every maximal ideal $P$ of $D$ is a $t$-ideal. Now consider the prime ideal $%
P[X]$ in $R[X]$ and note that $P[X]$ can never be a maximal ideal because $%
R[X]/P[X]\cong (R/P)[X]$ is a polynomial ring over a field and so must have
an infinite number of maximal ideals. This forces $P[X]$ to be properly
contained in an infinite number of maximal ideals $M_{\alpha }$ of $R[X].$
Let $M$ be one of them$.$ Then $M=(f,P[X]).$ Now, if it were the case that $%
R=(R,\mathcal{L}_{d}(R)\backslash \{0\},\Pi _{R})$ for $\Pi _{R}$ the set of
proper $t$-ideals, then every maximal ideal of $R$ would be a $t$-ideal.
This would make $M$ a $t$-ideal with $M\cap D=P\neq (0).$ But then,
according to Proposition 1.1 of \cite{HZ 1989}, $M=(M\cap D)[X]=P[X],$ a
contradiction to the fact that $P[X]\subsetneq M$. For (2) note that if $%
D=(D,\mathcal{L}_{t}(D),\Pi _{D})$ for $\Pi _{D}$ the set of proper nonzero
principal ideals (resp., $t$-invertible $t$-ideals, $t$-ideals of finite
type, divisorial ideals), then in each case every maximal $t$-ideal $P$ of $%
D $ is divisorial. Now let $M$ be a maximal $t$-ideal of $R.$ If $M\cap
D=(0) $, then $M$ is a $t$-invertible $t$-ideal and hence divisorial by
Theorem 1.4 of \cite{HZ 1989}. Next if $M$ is such that $M\cap D\neq (0),$
then $M=(M\cap D)[X]$ where $M\cap D$ is a maximal $t$-ideal of $D$ and
hence divisorial. Now it is easy to show that $M$ is divisorial. Conversely
suppose that $R=(R,\mathcal{L}_{t}(R),\Pi _{R})$ for $\Pi _{R}$ the set of
proper divisorial ideals. Then every maximal $t$-ideal $M$ of $R$ is
divisorial. Now let $P$ be a maximal $t$-ideal of $D.$ Then $P[X]$ is a
maximal $t$-ideal of $R$ by Proposition 1.1 of \cite{HZ 1989} and hence
divisorial. But this leads to $P[X]=(P[X])_{v}=P_{v}[X]$ and hence to $%
P=P_{v}.$
\end{proof}

Before checking how fare the $D=(D,\mathcal{L}_{d}(D)\backslash \{0\},\Pi
_{D})$ for $\Pi _{D}$ the set of proper nonzero principal ideals (resp., $t$%
-invertible $t$-ideals, $t$-ideals, $t$-ideals of finite type divisorial
ideals) under extension to rings of fractions, we need to consider another
type of extension, the $D+XL[X]$ construction. Yet to be able to fully
appreciate how it works, one needs to learn a little about it, the
construction $D+XL[X]$. Let $D$ be an integral domain with quotient field $K$%
, let $L$ be an extension of $K$ and let $X$ be an indeterminate over $L$.
Then $R=D+XL[X]$ $=\{f\in L[X]|f(0)\in D\}$ is an integral domain.

Indeed $R$ has two kinds of prime ideals $P$, ones that intersect $D$
trivially and ones that don't. If $P\cap D\neq (0)$ then $P=P\cap D+XL[X]$ 
\cite[Lemma 1.1]{CMZ 1986} and obviously $P$ is maximal if and only if $%
P\cap D$ is. It can be shown, as was indicated prior to the proof of
Corollary 16 in \cite{ACZ 2015}, that if $P=P\cap D+XL[X]$, then $P$ is a
maximal $t$-ideal of $R$ if and only if $P\cap D$ is a maximal $t$-ideal of $%
D$ and indeed as $P_{v}=(P\cap D)_{v}+XL[X]$, $P$ is divisorial if and only
if $(P\cap D)$ is. Moreover, prime ideals of $R$ that are not comparable to $%
XL[X]$ are of the form $(1+Xg(X))R$ where $1+Xg(X)$ is an irreducible
element of $L[X]$, \cite[Lemmas 1.2, 1.5]{CMZ 1986}. Also as $XL[X]$ is of
height one $XL[X]$ is a $t$-ideal.

\begin{proposition}
\label{Proposition UXC}

\begin{enumerate}
\item Let $L$ be an extension of the field of quotients $K$ of an integral
domain $D$. let $X$ be an indeterminate over $L$ and let $R=D+XL[X]$. Then $%
D=(D,\mathcal{L}_{d}(D)\backslash \{0\},\Pi _{D})$ for $\Pi _{D}$ the set of
proper nonzero principal ideals (resp., $t$-invertible $t$-ideals, $t$%
-ideals, $t$-ideals of finite type divisorial ideals) if and only if $R=(R,%
\mathcal{L}_{d}(R)\backslash \{0\},\Pi _{R})$ with $\Pi _{R}$ the set of
proper nonzero principal ideals (resp., $t$-invertible $t$-ideals, $t$%
-ideals, $t$-ideals of finite type divisorial ideals) of $R$.

\item Let $L$ be an extension of the field of quotients $K$ of an integral
domain $D$. let $X$ be an indeterminate over $L$ and let $R=D+XL[X]$. Then $%
D=(D,\mathcal{L}_{t}(D),\Pi _{D})$ for $\Pi _{D}$ the set of proper nonzero
principal ideals (resp., $t$-invertible $t$-ideals, $t$-ideals of finite
type, divisorial ideals) if and only if $R=(R,\mathcal{L}_{t}(R),\Pi _{R})$
with $\Pi _{R}$ the set of proper nonzero principal ideals (resp., $t$%
-invertible $t$-ideals, $t$-ideals of finite type divisorial ideals) of $R$.
\end{enumerate}
\end{proposition}

\begin{proof}
(1) We only carry the proof through for one case of $\Pi _{D}$ leaving the
rest to the reader. Let $D=(D,\mathcal{L}_{d}(D)\backslash \{0\},\Pi _{D})$
for $\Pi _{D}$ the set of proper $t$-ideals (resp. $t$-ideals of finite
type). Then by Proposition \ref{Proposition UX1}, every maximal ideal of $D$
is a $t$-ideal (resp. $t$-ideal of finite type) and to show that $R=(R,%
\mathcal{L}_{d}(R)\backslash \{0\},\Pi _{R})$ with $\Pi _{R}$ the set of
proper $t$-ideals (resp. $t$-ideals of finite type) all we need do is check
if every maximal ideal $M$ of $R$ is a $t$-ideal (resp. $t$-ideal of finite
type). Now, according to the discussion prior to this proposition, any
maximal ideal incomparable to $XL[X]$ is principal and hence a $t$-ideal (of
finite type) and any maximal ideal $M$ comparable to $XL[X]$ is of the form $%
M=P+XL[X]$ where $P$ is a maximal ideal if and only if $M$ is and $P$ is a $%
t $-ideal if and only if $M$ is. Since every maximal ideal of $D$ is a $t$%
-ideal (resp. $t$-ideal of finite type) $P$ is a $t$-ideal (resp. $t$-ideal
of finite type) and it is easy to see that so is $M$ a $t$-ideal (resp. $t$%
-ideal of finite type). For the converse note that for every maximal ideal $%
P $ of $D$, $M=P+XL[X]$ is a maximal ideal of $R$ and so is a $t$-ideal
(resp. $t$-ideal of finite type), but this forces $P$ to be a $t$-ideal
(resp. $t$-ideal of finite type). We shall do one part of the proof of (2)
and leave the rest to the reader. Suppose that $D=(D,\mathcal{L}_{t}(D),\Pi
_{D})$ for $\Pi _{D}$ the set of proper divisorial ideals of $D.$ Then every
maximal $t$-ideal $P$ is a divisorial ideal of $D.$ Now every maximal $t$%
-ideal of $R$ that is incomparable with $XL[X]$ is a prime ideal, hence a
principal maximal ideal of the form $(1+Xg(X))R$ and hence a divisorial
ideal. This leaves maximal $t$-ideal of $R$ that are comparable with $XL[X].$
These are of the form $M=P+XL[X]$ where $P$ is a maximal $t$-ideal of $D.$
But maximal $t$-ideals of $D$ are divisorial and $%
M_{v}=(P+XL[X])_{v}=P_{v}+XL[X]=P+XL[X] $ (by discussion prior to this
proposition). Thus $R=(R,\mathcal{L}_{t}(R),\Pi _{R})$ with $\Pi _{R}$ the
set of proper nonzero divisorial ideals of $R.$ For the converse let $P$ be
a maximal $t$-ideal of $D$. Then $M=P+XL[X]$ is a maximal $t$-ideal. But, by
the condition, $M$ is divisorial which forces $P$ to be divisorial. Whence $%
D=(D,\mathcal{L}_{t}(D),\Pi _{D})$ for $\Pi _{D}$ the set of proper
divisorial ideals of $D.$
\end{proof}

Now we are ready to show that if $R=D_{S}$, for a multiplicative set $S$ of $%
D=(D,\mathcal{L}_{d}(D)\backslash \{0\},\Pi _{D})$ for $\Pi _{D}$ the set of
proper nonzero principal ideals (resp., $t$-invertible $t$-ideals, $t$%
-ideals, $t$-ideals of finite type, divisorial ideals), then it may not
generally be the case that $R=(R,\mathcal{L}_{d}(R)\backslash \{0\},\Pi
_{R}) $ for $\Pi _{R}$ the set of proper nonzero principal ideals (resp., $t$%
-invertible $t$-ideals, $t$-ideals, $t$-ideals of finite type, divisorial
ideals) of $R$. Let's first recall from Lemma \ref{Lemma TX2} that if every
maximal ideal is a $t$-invertible $t$-ideal then every maximal ideal is
actually invertible. Now let's start constructing examples.

\begin{example}
\label{Example UXD}Let $L$ be field extension of $K$ with $[L:K]=\infty $,
let $X$ be an indeterminate over $L$ and consider $R=D+XL[X]$. Set $%
S=D\backslash \{0\}$. If every maximal ideal of $D$ is principal
(invertible, $t$-ideal of finite type then so is every maximal ideal of $R$.
But that is not the case for every maximal ideal of $R_{S}$. For $%
R_{S}=K+XL[X]$ has a maximal ideal that is a $t$-ideal but not of finite
type and hence not invertible, nor principal.
\end{example}

The following example has been taken, almost verbatim, from \cite[Example 3.3%
]{HZ 2015}.

\begin{example}
\label{Example UXE} There does exist at least one direct example of a domain 
$D$ such that each maximal ideal of $D$ is a $t$-ideal but for some maximal $%
M$ we have $MD_{M}$ not a $t$-ideal. One such example is that of an
essential domain that is not a PVMD. (Recall that an integral domain $D$ is
essential if $D$ has a set $F$ of primes such that $D_{p}$ is a valuation
domain for each $P\in F$ and $D=\cap _{P\in F}D_{P}.)$. Now the example in
question was constructed by Heinzer and Ohm in \cite{HO 1973} and further
analyzed in \cite{MZ 1981} and \cite{GBL 2004}. As it stands the example has
all except one maximal ideals of height one primes and hence $t$-ideals and
the other maximal ideal $M$ is a height $2$ prime $t$-ideal. Indeed this is
the maximal ideal $M$ such that $D_{M}$ is a $2$-dimensional regular local
ring and so with a maximal ideal that is not a $t$-ideal. Showing that while 
$D=(D,\mathcal{L}_{d}(D)\backslash \{0\},\Pi _{D})$ for $\Pi _{D}$ the set
of $t$-ideals of $D$, $D_{M}\neq (D_{M},\mathcal{L}_{d}(D_{M})\backslash
\{0\},\Pi _{D_{M}})$ for $\Pi _{D_{M}}$ the set of $t$-ideals of $D_{M}$.
\end{example}

Now the fact that $D=(D,\mathcal{L}_{t}(D)$,$\Pi _{D})$ can go through the $%
D+XL[X]$ construction with the various definitions of $\Pi _{D}$ can be used
to construct, for example a domain of any dimension with $t$-maximal ideals
principal. If that reminds an attentive reader of comments (3) and (4) of
Remarks 8 of \cite{MZ 1990}, then so be it. The point however is that the
domains 
\begin{equation*}
D=(D,\mathcal{L}_{d}(D)\backslash \{0\},\Pi _{D})
\end{equation*}
and $D=(D,\mathcal{L}_{t}(D),\Pi _{D}),$ with suitable $\Pi _{D}$s, do not
have the usual Ascending Chain Conditions on ideals (principal or $t$%
-)ideals. One may wonder if there are any simple restrictions that will get
the beast under control. Yet to prepare to see that, here is another simple
set of results that can come in handy when we are dealing with completely
integrally closed integral domains.

Of course before we bring in those results some introduction is in order.
Recall that an integral domain $D$ with quotient field $K$ is completely
integrally closed if whenever $rx^{n}$ $\in D$ for $x\in K$, $0\neq $ $r\in D
$, and every integer $n\geq 1$, we have $x\in D$. It can be shown that an
intersection of completely integrally closed domains is completely
integrally closed. The go to reference for Krull domains is Fossum's book 
\cite{Fos 1973} where you can find that $D$ is a Krull domain if $D$ is a
locally finite intersection of localizations at height one primes such that $%
D_{P}$ is a discrete valuation domain at each height one prime. Thus a Krull
domain is completely integrally closed. 

Glaz and Vasconcelos \cite{GV 1977} called an integral domain $D$ an
H-domain if there is an ideal $A$ with $A^{-1}=D$, (or equivalently $A_{v}=D)
$ then $A$ contains a finitely generated subideal $F$ such that $%
A^{-1}=F^{-1}$. They showed that a completely integrally closed H-domain is
a Krull domain. In \cite[Proposition 2.4]{HZ 1988} it was shown that $D$ is
an H-domain if and only if every maximal $t$-ideal of $D$ is divisorial. We
have in the following a basic result and some of its derivatives.

\begin{proposition}
\label{Proposition UX1}Let $D$ be a completely integrally closed domain. Then

\begin{enumerate}
\item $D=$ $(D,\mathcal{L}_{t}(D),\Pi )$ is a Krull domain if and only if $%
\Pi $ is the set of proper divisorial ideals of $D$.

\item $D=(D,\mathcal{L}_{t}(D),\Pi )$ is a locally factorial Krull domain if
and only if $\Pi $ is the set of proper invertible integral ideals of $D$.

\item $D=(D,\mathcal{L}_{t}(D),\Pi )$ is a Krull domain if and only if $\Pi $
is the set of proper $t$-invertible $t$-ideals of $D$.

\item Let $D$ be such that $D_{M}$ is a Krull domain for each maximal ideal $%
M$ of $D$. Then $D=$ $(D,\mathcal{L}_{t}(D),\Pi )$ is a Krull domain if and
only if $\Pi $ is the set of proper divisorial ideals of $D$ \cite{EIT 2019}.

\item Let $D$ be an intersection of rank one valuation domains. Then $D=$ $%
(D,\mathcal{L}_{t}(D),\Pi )$ is a Krull domain if and only if $\Pi $ is the
set of proper divisorial ideals of $D$.

\item Let $D$ be an almost Dedekind domain. Then $D=$ $(D,\mathcal{L}%
_{d}(D),\Pi )$ is a Dedekind domain if and only if $\Pi $ is the set of
proper divisorial ideals of $D$.
\end{enumerate}
\end{proposition}

\begin{proof}
The idea of proof, in each case, is that every maximal $t$-ideal (maximal
ideal) being contained in a proper divisorial ideal must be equal to it and
combining this with the fact that $D$ is completely integrally closed we get
the Krull domain conclusion. For the locally factorial domain conclusion in
(2) we note that every maximal $t$-ideal of $D$ is invertible and so
divisorial. This gives the Krull conclusion and a Krull domain is locally
factorial if and only if every height one prime of $D$ is invertible \cite[%
Theorem 1]{And 1978}. For the Dedekind domain conclusion in (6), we note
that every maximal ideal is of height one and divisorial, being invertible.
So every maximal ideal is a $t$-ideal and so the domain is Krull and one
dimensional. The converse in each case is obvious, in that if $D$ is a Krull
domain then $D$ is completely integrally closed and every maximal $t$-ideal
of $D$ is, a $t$-invertible $t$-ideal and hence, divisorial. (If $D$ is
locally factorial, as in (2), every maximal $t$-ideal of $D$ is invertible
and hence divisorial.) And if $D$ is Dedekind, then $D$ is completely
integrally closed and every maximal ideal is invertible and hence divisorial.
\end{proof}

\begin{theorem}
\label{Theorem UX2}Let $D$ be a Noetherian domain (resp., has ACC on
invertible ideals). Then

\begin{enumerate}
\item $D=$ $(D,\mathcal{L}_{d}(D)\backslash \{0\},\Pi )$ is a PID if and
only if $\Pi $ is the set of proper nonzero principal ideals of $D$.

\item $D=(D,\mathcal{L}_{d}(D)\backslash \{0\},\Pi )$ is a Dedekind domain
if and only if $\Pi $ is the set of proper invertible ideals of $D$.
\end{enumerate}
\end{theorem}

\begin{proof}
$D=$ $(D,\mathcal{L}_{d}(D)\backslash \{0\},\Pi )\Leftrightarrow \forall
A\in \mathcal{L}_{d}(D)\backslash \{0\}$ $(A\neq D$ $\Rightarrow \exists \pi
\in \Pi (A\subseteq \pi ))$ where $\Pi $ is the set of proper nonzero
principal (invertible) ideals. Start with a proper nonzero ideal $A$ of $D.$
Then by the condition $A\subseteq \pi _{1}$ for some $\pi _{1}\in \Pi .$ Let 
$A_{1}=A\pi _{1}^{-1}.$ Then $A\subsetneq A_{1}.$ If $A_{1}=D$ we have $%
A=\pi _{1}$ a principal (invertible) ideal and we are done. If $A_{1}\neq D$
then by the condition $A_{1}\subseteq \pi _{2}$ where $\pi _{2}$ is a proper
principal (invertible) ideal. Let $A_{2}=A_{1}\pi _{2}^{-1}$ $=A(\pi _{1}\pi
_{2})^{-1}.$ Then $A\subsetneq A_{1}\subsetneq A_{2}$ and at a general stage 
$A\subsetneq A_{1}\subsetneq A_{2}\subsetneq ...\subsetneq A_{r}\subsetneq
...$ where $A_{i}=A(\Pi _{j=1}^{i}\pi _{j})^{-1}.$ Now, being a Noetherian
domain (or with ACC on invertible ideals), $D$ cannot afford an infinite
strictly ascending chain of proper invertible integral ideals. Whence at
some stage $n,$ $A_{n}=A(\Pi _{j=1}^{n}\pi _{j})^{-1}=D,$ forcing $A=\Pi
_{j=1}^{n}\pi _{j}$. This makes a typical ideal $A$ of $D$ principal
(invertible). The converse in both cases is obvious.
\end{proof}

The ideas touched on above got developed in a different direction, in \cite%
{Zaf 2020}, using predicates; a demonstration of the fact that the change of
language is at times the change of degrees of freedom, in some direction.
The point of this exercise here is to indicate how pre-Riesz monoids can
help create examples. Granted that some of the examples are a repeat of what
we already know, but we cannot kill a new idea just because the initial
examples are gathered from known sources. Besides when you are in a
demonstration mode, known examples are actually welcome and useful in fixing
the ideas. My response to any, possibly, raised eyebrows is that there are
examples of direct translations from ring theoretic results to results on
monoids in Halter-Koch's book \cite{FHK B} and in most of his work close to
the end of his career. (See for instance his paper on mixed invertibility 
\cite{FHK mix}.) There are other examples, even in ideal theory, almost all
results on PVMDs, originally proved using the $t$-operation are now being
stated and proved for the so-called $w$-operation. So, hopefully, there's no
problem with some examples being duplicated.

There is also the thought that as \textquotedblleft contains means
divides\textquotedblright\ some of these ideas may be extended to
non-commutative systems and/or to nearrings. Also the \textquotedblleft
designation\textquotedblright\ of $\Pi _{D}$ maybe helpful in introducing an
element of randomness, such as \textquotedblleft let's see if $\Pi _{D}$
consists of only two ideals\textquotedblright\ or if \textquotedblleft every
ideal is contained in at least one and at most a finite number of members of 
$\Pi _{D}$\textquotedblright . We can also take a set of ideals $S$ and
designate $\Pi _{D}$ as the set of ideals that are finite products of
members of $S$. Finally, if the notation irks some readers, I adopted the
notation $D=(D,\mathcal{L}_{t}(D),\Pi )$ etc., because it enables me to
state and prove a number of theorems in one go. It is also a shorthand way
of recording one's thoughts about domains of this type. If a reader has some
other ideas about such a device, the reader is invited to go ahead and
produce a more nimble shorthand.

\section{Virtual factoriality \label{S4}}

Now let me end the article with a mention of something positive, a kind of
unique factorization. Recall that two positive elements $a,b$ of a p.o.
group $G$ are said to be disjoint if $h\leq a,b$ implies $h\leq e$, the
identity of the group. While this definition works admirably in p.o. groups,
it can cause confusion in the monoid set up. If, following Birkhoff, (see
comment on page 220 of \cite{B 1948}),we liken disjointness with
\textquotedblleft relative primeness\textquotedblright\ or coprimeness, then
in an integral domain there are two kinds of coprimeness: One that comes
from the lack of non-unit common factors in a set of elements and one that
comes from the elements being disjoint in the group of divisibility of the
domain. As the theory of factorization in integral domains developed, the
two kinds got distinct names.

\begin{definition}
\label{Definition XX0}Let $x,y$ be two elements in $D\backslash \{0\}$,
where $D$ is an integral domain. Then $x,y$ are coprime if $z|x,y$ implies
that $z$ is a unit in $D$, i.e., if $GCD(x,y)$ is a unit and $x,y$ are $v$%
-coprime if $xD\cap yD=xyD,$ i.e., if $(x,y)_{v}=D$.
\end{definition}

As indicated in \cite{Zaf 2006}, page 389, the notion of $x,y$ being $v$%
-coprime gets translated to $xD,yD$ being disjoint in $G(D)=\{xD|x\in
K\backslash \{0\}\}$ ordered by $xD\leq yD$ $\Leftrightarrow xD\supseteq yD$%
. As is also indicated in \cite{Zaf 2006} there are several shades of
coprimality, thanks to the star operations. We may say that $x,y$ are $\star 
$-coprime if $(x,y)^{\star }=D$. But the kinds of coprimality outlined in
the definition above are the two extremes. Again, it was indicated in \cite%
{Zaf 2006} that while $x,y\in D\backslash \{0\}$ being $v$-coprime implies
that $x,y$ are coprime, it is not the case the other way round. That is the
notion of $v$-coprimality is more reliable. And in a general monoid
situation our best bet would be to stick with the most reliable, unless a
specific situation requires a deeper digging. It turns out that the best use
of $v$-coprimality comes from the results that if, for $a$,$b$,$c\in
D\backslash \{0\},$ $a$ and $c$ are $v$-coprime and if $a|bc$, then $a|b$,
and $(a,bc)_{v}=D$ if and only if $(a,b)_{v}=(a,c)_{v}=D$, see Propositions
2.2 and 2.3 of \cite{Zaf 2006}. We shall retain the definition for
disjointness used for p.o. groups, that is we call two positive elements $a,b
$ in a p.o. monoid $M$ disjoint if for every $c\leq a,b$ we have $c\leq e$,
the identity of $M$, i.e., $a\wedge b=0$. We may apply adjustments where
necessary.

We note that if in a Riesz monoid $\langle M,\ast ,e,\leq \rangle $ we have $%
e\leq a\leq b\ast c$ where $a\wedge c=e$, i.e., $a$ and $c$ are disjoint,
then $a\leq b$. This can be shown by noting that $a$ is primal. Generally if
for $\langle M,\ast ,e,\leq \rangle $, $M^{+}$ is the positive cone of a
directed p.o. group $G$ one can prove the same result, because if $a\wedge c$
(resp., $a\vee c)$ exists in $G$ then $x\ast (a\wedge c)=(x\ast a)\wedge
(x\ast c)$ (resp., $x\ast a\vee c=(x\ast a)\vee (x\ast c))$ for all $x\in G$
and hence for all $x\in G^{+}$, see e.g. \cite{YZ 2011}.

\begin{proposition}
\label{Proposition XX1}Let $M^{+}$ be the positive cone of a directed p.o.
group that is abelian. Then in $M^{+}\backslash \{e\}$ $a\leq b\ast c$,
where $a\wedge c=e$, in $G$, implies that $a\leq b$.
\end{proposition}

\begin{proof}
Note that as $a\wedge c=e,$ in $G,$ we have $b=b\ast (a\wedge c)=(b\ast
a)\wedge (b\ast c),$ by the remark prior to the statement of this
proposition. Now $a\leq a\ast b,b\ast c$ implies $a\leq (b\ast a)\wedge
(b\ast c)=b\ast (a\wedge c)=b\ast e=b.$
\end{proof}

A result similar to the above can also be proved in, what we have chosen to
call, the reverse multiplicative lattice $\Delta _{\star }(D)=\langle
I^{\star }(D),+^{\star },\times ^{\star },\wedge ,\vee ,D\rangle $.

\begin{proposition}
\label{Proposition XXII}If $A,B,C\in $ $\Delta _{\star }(D)=\langle I^{\star
}(D),+^{\star },\times ^{\star },\wedge ,\vee ,D\rangle $ such that $A\wedge
C=D$, then $A\leq B\times ^{\star }C$ implies that $A\leq B$.
\end{proposition}

\begin{proof}
Recall that reverse containment rules in $\Delta _{\star }(D)$. Consequently 
$D$ is the least element of the lattice. Thus $X\leq Y$ in $\Delta _{\star
}(D)$ if and only if $X\supseteq Y$ and $X\wedge Y$ becomes $(X+^{\star
}Y)=(X+Y)^{\star }.$ Thus the statement "$A\leq B\times ^{\star }C$ where $%
A\wedge C=D"$ translates to "$A\supseteq B\times ^{\star }C$ where $%
(A+C)^{\star }=D".$ Now as we also have $A\supseteq A\times ^{\star }B$ we
can say that $A\supseteq A\times ^{\star }B,B\times ^{\star }C$ and so, as $%
A $ is a $\star $-ideal we have $A\supseteq A\times ^{\star }B+^{\star
}B\times ^{\star }C=(AB+BC)^{\star }$ $=(B(A+C))^{\star }=(B(A+C)^{\star
})^{\star }=B,$ giving us $A\supseteq B$ which translates back to $A\leq B.$
\end{proof}

Proposition \ref{Proposition XXII} serves several purposes. First of all,
because there is no concept of invertibility or cancellation in $\Delta
_{\ast }(D)$, nor in $\mathcal{L}_{\ast }(D)$, it shows that results such as
\textquotedblleft $a\leq b\ast c$ and $a\wedge c=e$, implies that $a\leq b$%
\textquotedblright\ may hold in monoids that are neither divisibility nor
cancellative. It also shows that while in some monoids $a\wedge b=e$ may
make sense as the elements $a,b$ being disjoint, i.e. $h\leq a,b$ implies $%
h\leq e$, in others such as $\mathcal{L}_{\ast }(D)$ $a,b$ being disjoint
may well be equivalent to $a\vee b=1$ where $1$ is the largest element of
the lattice. Of course in $\mathcal{L}_{\ast }(D)$ $a\wedge b=0$, if and
only if $a=0$ or $b=0$. Then, as pointed out, after definition 12, in \cite%
{LYZ 2020}, there are monoids in which $a\wedge b_{i}=e$ for $i=1,..,n$ but $%
a\wedge (b_{1}\ast ...\ast b_{n})>e.$

One purpose of studying monoids is to see if some kind of unique
factorization exists in them and it helps to separate or collect factors of
an element of a monoid by disjointness or lack of it, in the monoid
environment. So if in a monoid $M$, for elements $a,b_{1},...,b_{n}$ such
that $a\wedge b_{i}=e$, for $i=1,...,n,$ implies $a\wedge (b_{1}\ast ...\ast
b_{n})=e$ we know that we are in a safe environment. This can happen in
cases indicated in the following proposition. Yet before that let's call $M$
a multiplicative $\wedge $-semilattice if $a\wedge b\in M$ for all $a,b\in M$
and for all $x\in M$ we have $x\ast (a\wedge b)=x\ast a\wedge x\ast b$, $%
\wedge $ is an associative and commutative binary operation and $x=x\wedge x$
for all $x\in M$. Also call a monoid $M$ $\wedge $-smooth if whenever $%
a\wedge b$ exists for some $a,b\in M$ we have $x\ast (a\wedge b)=x\ast
a\wedge x\ast b$ for all $x\in M$ and if $a\wedge (x\wedge y)$ or $(a\wedge
x)\wedge y)$ exists we have $a\wedge (x\wedge y)=(a\wedge x)\wedge y)$.

\begin{proposition}
\label{Proposition YX}

\begin{enumerate}
\item Let $M$ be a $\wedge $-smooth monoid. If, for $a,b_{1},...,b_{n}\in M$%
, we have $a\wedge b_{i}=e$, for $i=1,...,n$ in $M$, then $a\wedge
(b_{1}\ast ...\ast b_{n})=e$. Conversely if, for $a,b_{1},...,b_{n}\in M,$
we have $a\wedge (b_{1}\ast ...\ast b_{n})=e$, then $a\wedge b_{i}=e$ for
all $i=1,2,...,n$.

\item If $M$ is a multiplicative $\wedge $-semilattice and if, for $%
a,b_{1},...,b_{n}\in M,$ we have $a\wedge b_{i}=e$, for $i=1,...,n,$ in $M,$
then $a\wedge (b_{1}\ast ...\ast b_{n})=e$.
\end{enumerate}
\end{proposition}

\begin{proof}
(1). The case of $n\leq 2$ is provided by the definition. Suppose that we
have established the proposition for $n\leq r$ for $r\geq 2.$ That is if $%
a\wedge b_{i}=e$, for $i=1,...,r$ in $M,$ then $a\wedge (b_{1}\ast ...\ast
b_{r})=e.$ Let $a\wedge b_{r+1}=e$ and consider $a\wedge (b_{1}\ast ...\ast
b_{r+1}).$ Since $a\leq a\ast b_{1}\ast ...\ast b_{r}$ we have $a=a\wedge
a\ast b_{1}\ast ...\ast b_{r}$ and so $a\wedge (b_{1}\ast ...\ast
b_{r+1})=(a\wedge (a\ast b_{1}\ast ...\ast b_{r}))\wedge (b_{1}\ast ...\ast
b_{r+1})=$ $(a\wedge \left( (a\ast b_{1}\ast ...\ast b_{r}))\wedge
(b_{1}\ast ...\ast b_{r+1}\right) $ $=$ $a\wedge ((b_{1}\ast ...\ast
b_{r})\ast (a\wedge b_{r+1}))$ (and because $(a\wedge b_{r+1})=e$ we have)

$a\wedge (b_{1}\ast ...\ast b_{r+1})=a\wedge ((b_{1}\ast ...\ast b_{r}).$
But by the induction hypothesis, we have $a\wedge ((b_{1}\ast ...\ast
b_{r})=e,$ So $a\wedge (b_{1}\ast ...\ast b_{r+1})=e.$ Thus (1) holds for
all integers $>1.$ For the converse note that as $b_{i}\leq b_{1}\ast
b_{2}\ast ...\ast b_{n}$ we have $b_{i}=b_{i}\wedge (b_{1}\ast b_{2}\ast
...\ast b_{n}).$ Thus $a\wedge b_{i}=a\wedge (b_{i}\wedge (b_{1}\ast
b_{2}\ast ...\ast b_{n})$ $=a\wedge (b_{1}\ast b_{2}\ast ...\ast b_{n}\wedge
b_{i})$ $=(a\wedge b_{1}\ast b_{2}\ast ...\ast b_{n})\wedge b_{i}=e\wedge
b_{i}.$

(2). Follows from (1), because a multiplicative $\wedge $-semilattice is $%
\wedge $-smooth.
\end{proof}

Dually we can talk about multiplicative $\vee $-semilattices, and $\vee $%
-smooth monoids, and noting that the dual of $a\wedge b=e$ is $a\vee b=1$,
the largest element, we can state and prove the dual of Proposition \ref%
{Proposition YX}, by replacing $\wedge $ by $\vee $, $\leq $ by $\geq $ and $%
e$ by $1$.

Next, there are monoids such as the monoid of nonzero finitely generated
integral ideals of $D$, partially ordered by inclusion (or by reverse
inclusion) denoted by $f(D)$ (respectively by $\varphi (D))$ under
multiplication of ideals. We can also talk about $f_{\ast }(D)$ and $\varphi
_{\ast }(D)$ of integral $\ast $-ideals of finite type closed under $\ast $%
-multiplication.

Now in a pre-Riesz monoid, as we have defined it, there seems to be no
indication if a pre-Riesz monoid is $\wedge $-smooth (or $\vee $-smooth).
One way of dealing with the situation is to require that we deal only with
smooth Pre-Riesz monoids. Thus we have the obvious corollary. At the same
time it seems best to ask: Is it true that a pre-Riesz monoid is $\wedge $%
-smooth?

\begin{corollary}
\label{Corollary YX1}Let $M$ be a $\wedge $-smooth pre-Riesz monoid. If, for 
$a,b_{1},...,b_{n}\in M^{+}$, we have $a\wedge b_{i}=e$, for $i=1,...,n,$ in 
$M$, then $a\wedge (b_{1}\ast ...\ast b_{n})=e$. Conversely if $a\wedge
(b_{1}\ast ...\ast b_{n})=e,$ $b_{1},...,b_{n}\in M^{+}$, then $a\wedge
b_{i}=e$.
\end{corollary}

My interest in $\wedge $ (or $\vee )$-smooth pre-Riesz monoids arose from
the fact that they are as amenable to factorization as Riesz monoids. That
was the reason why I looked into the bases of pre-Riesz groups via Conrad's
F-condition, with Y.C. Yang \cite{YZ 2011}. I was hoping to find the
ultimate building blocks of factorization, in the positive cones of
pre-Riesz groups, as Conrad did in the form of a basic element in the case
of lattice ordered groups. Let's call an element $h$ of a monoid $M$ a
homogeneous element if $h>e$ (i.e. $h$ is strictly positive) and for all $%
u,v\in (e,h]$ (= for all $e<u,v\leq h)$ we have $e<l\leq u,v)$. It is easy
to see that if $x$ is a homogeneous element of $M$, then so is each $e<t\leq
x$. Let's first establish that homogeneous elements are not too hard to find.

\begin{lemma}
\label{Lemma YX2}Let $M$ be a pre-Riesz monoid such that $M^{+}=M$. Then the
following hold.

\begin{enumerate}
\item For $x,y,h\in M$, $e<h\leq x$ and $x\wedge y=e$ implies $h\wedge y=e$.

\item For $x,y,h_{1},h_{2}\in M$, $e<h_{1},h_{2}\leq x$ and $x\wedge y=e$
implies $h_{i}\wedge y=e$ and if in addition $h_{1}\wedge h_{2}=e$, then $%
h_{1},h_{2},y$ are mutually disjoint.

\item If $x_{1},...,x_{n}\in M$ such that $x_{i}$ are mutually disjoint and
we have $e<h_{i1},h_{i2}\leq x_{i}$ , for some $i\in \{1,...,n\}$, such that 
$h_{i1}\wedge h_{i2}=e$, then $%
\{x_{1},...,x_{i-1},h_{i1},h_{i2},x_{i+1},...,x_{n}\}$ are mutually disjoint.

\item If $e<h,k\in M$, then $h$ and $k$ are nondisjoint if and only if there
is $0<t\leq h,k$.
\end{enumerate}
\end{lemma}

\begin{proof}
(1) If $h,y$ are non-disjoint then by the pre-Riesz property, there is $%
e<r\leq h,y.$ But since $h\leq x$ we have $e<r\leq x,y$ a contradiction. Now
(2) follows from (1) by noting that each of $h_{i}$ is disjoint with $y$ and
if we add to it the fact that $h_{1}\wedge h_{2}=e$ we have the conclusion
that $\{h_{1},h_{2},y\}$ are mutually disjoint. For (3) note that $%
\{h_{i1,}h_{i2},x_{j}|i\neq j\}$ are mutually disjoint by (2). (4) follows
from the definition of a pre-Riesz monoid.
\end{proof}

\begin{proposition}
\label{Proposition YX3}Let $M$ be a pre-Riesz monoid such that $M^{+}=M$,
i.e. $M$ is conic. If $M$ satisfies \textquotedblleft CFC: every strictly
positive element exceeds at most a finite number of mutually disjoint
strictly positive elements\textquotedblright , then every strictly positive
element of $M$ exceeds at least one homogeneous element and at most a finite
number of mutually disjoint homogeneous elements.
\end{proposition}

\begin{proof}
Let $e<x\in M$ and let there be $e<u,v\leq x.$ The pre-Riesz condition
provides that if $u,v$ are non-disjoint there is a strictly positive $l$
such that $e<l\leq u,v\leq x.$ Now if $x$ is such that for all $u,v\in (e,h]$
we have $e<l\leq u,v,$ then $x$ is a homogeneous element. If $x$ is not
homogeneous then there must be at least two disjoint elements $u_{1},v_{1}$
preceding $x.$ Let $n$ be the largest number of mutually disjoint strictly
positive elements of $M$ preceding $x$ and let $\{x_{1},...,x_{n}\}$ be a
set of $n$ elements such that $0<x_{i}\leq x$ where $x_{i}$ are mutually
disjoint. Then each of $x_{i}$ is a homogeneous element. For if not, and say 
$x_{i}$ is not homogeneous, then there exist at least two elements $%
e<h_{i1},h_{i2}\leq x_{i}$ such that $h_{i1}\wedge h_{i2}=e.$ By Lemma \ref%
{Lemma YX2} the set $\{x_{1},...x_{i-1},h_{i1},h_{i2},x_{i+1},...,x_{n}\}$
of $n+1$ elements, consists of mutually disjoint strictly positive elements
preceding $x,$ a contradiction. Whence each of $x_{i}$ is homogeneous. Now
to establish beyond doubt that $x$ exceeds at least one homogeneous element
we proceed as follows. Let $e<x\in M$ and let there be $e<u,v\leq x.$ The
pre-Riesz condition provides that if $u,v$ are non-disjoint there is a
strictly positive $l$ such that $e<l\leq u,v\leq x.$ Now, as noted above, if 
$x$ is such that for all $u,v\in (e,h]$ we have $e<l\leq u,v,$ then $x$ is a
homogeneous element and we are done. If $x$ is not homogeneous then there is
a pair of strictly positive disjoint elements $h,k$ preceding $x.$ If either
of $h,k$ is homogeneous we can stop. If not, we can find a pair $h_{1},h_{2}$
of strictly positive disjoint elements preceding $h,$ noting that $%
h_{1},h_{2},k$ are mutually disjoint by Lemma \ref{Lemma YX2}, in particular 
$h_{2},k$ are disjoint. If neither of $h_{1},$ $h_{2}$ is homogeneous, find $%
e<h_{11},h_{12}\leq h_{1}$ to get $e<h_{11},h_{12},h_{2},k$ mutually
disjoint preceding $x$, noting that in particular $\{h_{12},h_{2},k$ $\}$
are mutually disjoint. Assuming that neither of $h_{11},h_{12}$ is
homogeneous and repeating the previous step with $h_{11}$ to get $%
e<h_{111},h_{112}\leq h_{11}$ such that $h_{111}\wedge h_{12}=e;$ noting
that by Lemma \ref{Lemma YX2} $\{h_{111},h_{112},h_{12},h_{2},k$ $\}$ are
mutually disjoint and in particular $\{h_{112},h_{12},h_{2},k$ $\}$ are
mutually disjoint. Now this cannot go on indefinitely as each new step
increases the number of mutually disjoint strictly positive elements
preceding $x$ by one and makes the number of mutually disjoint elements
preceding $x$ tend to infinity, contradicting CFC.
\end{proof}

\begin{remark}
\label{Remark YX4} 
\textquotedblleft CFC\textquotedblright in the statement of Proposition \ref%
{Proposition YX3} stands for Conrad's F-Condition, as used in \cite{Con 1961}%
. The plan in the proof of the latter part of Proposition \ref{Proposition
YX3} is taken from Theorem 5.2 of \cite{Con 1961}.
\end{remark}

The notion of a homogeneous element was developed in \cite{MRZ 2008} for the
study of factoriality in general Riesz groups, without any reference to
commutativity. So it should work well in a commutative Riesz monoid. We now
proceed to show that the notion of factoriality based on homogeneous
elements in Riesz groups works well in some pre-Riesz monoids. As a first
step, let's call two homogeneous elements $h,k\in M$ similar, denoted as $%
h\sim k$, if there is $e<t\leq h,k,$ that is if $h$ and $k$ are non-disjoint 
$(h\wedge k\neq e)$. Indeed if we use $(e,h]$ to mean the set of all $%
e<u\leq h$ and if $h$ is a homogeneous element of $M$, then each $u\in (e,h]$
is homogeneous.

\begin{proposition}
\label{Proposition YX5}Let $h$ and $k$ be two homogeneous elements in a
pre-Riesz monoid $M$. Then the following are equivalent.

\begin{enumerate}
\item $h\wedge k=e$.

\item For every pair $(a,b)\in (e,h]\times (e,k]$ we have $a\wedge b=e$.

\item For some pair $(a,b)\in (e,h]\times (e,k]$ we have $a\wedge b$ $=e$.
\end{enumerate}
\end{proposition}

\begin{proof}
(1) $\Rightarrow $ (2) Suppose (2) does not hold. Then for some pair $%
(a,b)\in (e,h]\times (e,k]$ we have $a\wedge b\neq e$ which, in a pre-Riesz
monoid, means that there is $e<t\leq a,b.$ But then $e<t\leq h,k$ which
forces $h\wedge k\neq e.$

(2) $\Rightarrow $ (3). Obvious.

(3) $\Rightarrow $ (1). Suppose that $h\wedge k\neq e.$ Then there is $%
e<t\leq h,k.$ This makes $t$ a homogeneous element. Now let $(a,b)\in
(e,h]\times (e,k].$ Since $e<t,a\leq h$ and $h$ is homogeneous we have $%
t\wedge a\neq e.$ Also since we are in a pre-Riesz monoid, there is a $%
e<t_{1}\leq t,a.$ Similarly there is a $e<t_{2}\leq t,b.$ Now since $t$ is
homogeneous, $e<t_{1},t_{2}\leq t$ and we are in a pre-Riesz monoid, there
must be a $e<t_{3}\leq t_{1},t_{2}.$ But then $t_{3}\leq a,b,$ forcing $%
a\wedge b\neq e.$ In other words, negation of (1) implies the negation of
(3).
\end{proof}

By negating the constituent statements in Proposition \ref{Proposition YX5}
we get the following statement.

\begin{proposition}
\label{Proposition YX6} Let $h$ and $k$ be two homogeneous elements in a
pre-Riesz monoid $M$. Then the following are equivalent.

\begin{enumerate}
\item $h\wedge k\neq e$.

\item For every pair $(a,b)\in (e,h]\times (e,k]$ we have $a\wedge b\neq e$.

\item For some pair $(a,b)\in (e,h]\times (e,k]$ we have $a\wedge b\neq e$.
\end{enumerate}
\end{proposition}

\begin{proposition}
\label{Proposition YX7} Let $h(M)$ be the set of all homogeneous elements of
a pre-Riesz monoid $M$. Then similarity is an equivalence relation on $H$.
\end{proposition}

\begin{proof}
The proof follows exactly the same lines as the proof of (6) of Proposition
1.1 of \cite{MRZ 2008}. That is, reflexivity and symmetry being clear, all
we need is check if transitivity works. For this let $u,v,w\in h(M)$ such
that $u\sim v$ and $v\sim w.$ Now as $u\sim v$ we have $u\wedge v\neq e.$ By
the pre-Riesz property there is $e<t_{1}\leq u,v.$ Next as $v\sim w$ we have 
$v\wedge w\neq e$ and by the pre-Riesz property there is $e<t_{2}\leq v,w.$
Now as $e<t_{1},t_{2}\leq v$ and as $v$ is homogeneous, and $t_{1},t_{2}$
strictly positive, we have $t_{1}\wedge t_{2}\neq e.$ But the pre-Riesz
property again gives $e<t_{3}\leq t_{1},t_{2}$ which implies that $%
e<t_{3}\leq t_{1}\leq u$ and $e<t_{3}\leq t_{2}\leq w.$ Implying $%
e<t_{3}\leq u,w$ which is the same as saying that $u\wedge w\neq e$ or $%
u\sim w.$
\end{proof}

\begin{proposition}
\label{Proposition YX8} Let $M$ be a pre-Riesz monoid. Let $e<x,y\in M$ such
that $x\wedge y=e$. If $h$ is a homogeneous element, then $h$ must be
disjoint with at least one of $x,y$. Generally if $x_{1},...,x_{n}$ are
mutually disjoint strictly positive elements of a pre-Riesz monoid $M$ and $%
h $ is a homogeneous element of $M$, then $h$ must be disjoint with at least 
$n-1$ of $x_{i}$. Consequently if $x_{1},...,x_{n}$ are mutually disjoint
strictly positive elements of a $\wedge $-smooth pre-Riesz monoid $M$ and $h$
is a homogeneous element of $M$ with $h\leq x_{1}\ast ...\ast x_{n}$, then $%
h\leq x_{i}$ for exactly one $1\leq i\leq n$.
\end{proposition}

\begin{proof}
Suppose on the contrary that $h$ is non-disjoint with both of $x$ and $y$.
Then since we are working inside a pre-Riesz monoid there are $e<t_{1}\leq
x,h$ and $e<t_{2}\leq y,h.$ But then $e<t_{1},t_{2}\leq h$ and $h$ is
homogeneous. This leads to the existence of $e<t_{3}\leq t_{1},t_{2}$ and to 
$e<t_{3}\leq x,y$ which contradicts the disjointness of $x$ and $y.$ For the
general case let $x=x_{i}$ and $y=y_{j},i\neq j$ . As we are working in a $%
\wedge $-smooth environment where $h\leq b\ast c$ and $h\wedge c=e$ implies $%
h\leq b$ and as, being homogeneous $h$ cannot be non-disjoint with more than
one disjoint elements, we can, by assuming $x_{i}\wedge h=e$ for $i\neq j,$
conclude that $h\leq x_{j}$ for exactly one $j.$
\end{proof}

\begin{proposition}
\label{Proposition YX9}Let $h_{1},h_{2}$ be two similar homogeneous elements
in a $\wedge $-smooth pre-Riesz monoid $M$. Then $h_{1}\ast h_{2}$ is a
homogeneous element similar to both $h_{i}$. Generally if $%
h_{1},h_{2},...,h_{n}$ are mutually similar homogeneous elements of a $%
\wedge $-smooth pre-Riesz monoid $M$, then $h_{1}\ast ...\ast h_{n}$ is a
homogeneous element similar to each of $h_{i}$.
\end{proposition}

\begin{proof}
Let $e<u,v\leq h_{1}\ast h_{2}.$ Claim that $u\wedge h_{i}\neq e$ and $%
v\wedge h_{i}\neq e.$ For if, say $u\wedge h_{1}=e,$ then $u\leq h_{1}\ast
h_{2}$ implies that $u\leq h_{2}.$ But as $e<u\leq h_{2}$ and $h_{2}$ is
homogeneous with $h_{1}\wedge h_{2}\neq e$ we must have $u\wedge h_{1}\neq e$%
, by Proposition \ref{Proposition YX6}, a contradiction. Now $u\wedge
h_{1}\neq e$ implies the existence of $e<t_{1}\leq u,h_{1}$, $v\wedge
h_{1}\neq e$ implies the existence of $e<t_{2}\leq v,h_{1}$ and as $h_{1}$
is homogeneous, $t_{1}\wedge t_{2}\neq e.$ But this, by the pre-Riesz
property, means that there is $e<t_{3}\leq t_{1},t_{2}$. That $e<t_{3}\leq
u,v$ is now obvious. Thus for any pair of strictly positive elements $u,v$
preceding $h_{1}\ast h_{2}$ we have $u\wedge v\neq e$ and so $h_{1}\ast
h_{2} $ is homogeneous. Now suppose that we have established the general
statement for $n\leq r.$ Then, by the induction hypothesis, $H=h_{1}\ast
...\ast h_{r}$ is homogeneous similar to each of $h_{i}.$ If $h_{r+1}$ is
homogeneous similar to any of $h_{i},$ $h_{r+1}$ is similar to $H$, by
Proposition \ref{Proposition YX6}. By the case of $n=2$, $H\ast h_{r+1}$ is
homogeneous similar to $H$ and $h_{r+1}$ and hence to all of $h_{i}.$
\end{proof}

Having established all the requirements of a theory of factoriality we
proceed as follows.

\begin{theorem}
\label{Theorem YX10} Let $M$ be a $\wedge $-smooth pre-Riesz monoid and let $%
x\in M$. If $x$ is expressible as $h_{1}\ast ...\ast h_{r}$ a sum/product of
finitely many homogeneous elements, then $x$ is expressible, uniquely, as a
sum/product of mutually disjoint homogeneous elements, up to permutation of
summands/factors.
\end{theorem}

Let $x=h_{1}\ast ...\ast h_{r}$, where $h_{i}$ are homogeneous. Pick, say, $%
h_{1}$ and pick all the homogeneous factors/summands similar to $h_{1}$.
Suppose that, by a relabeling, the first $n_{1}$ factors/summands of $x$ are
similar to $h_{1}$, the rest are disjoint with $h_{1}$ because similarity is
an equivalence relation. Let $H_{1}=h_{1}\ast ...\ast h_{n_{1}}$. Then $%
x=H_{1}\ast h_{n_{1}+1}\ast ...\ast h_{r}$, where $h_{n_{1}+1},...,h_{r}$
are all disjoint with $h_{1}$ and hence with $H_{1}$, via Proposition \ref%
{Proposition YX5}. Now repeat the previous step with $h_{n_{1}+1}$,
collecting all the homogeneous summands/factors similar to $h_{n_{1}+1}$ and
assume, relabeling if necessary, that $h_{n_{1}+1},...,h_{n_{2}}$ is the set
of all homogeneous factors of $h_{n_{1}+1}\ast ...\ast h_{r}$, (and hence of 
$x)$ similar to $h_{n_{1}+1}$ and form $H_{2}=h_{n_{1}+1}\ast ...\ast
h_{n_{2}}$, via Proposition \ref{Proposition YX9}. Since each of $h_{i}$ $%
(n_{1}+1\leq i\leq n_{2})$ is disjoint with $h_{1}$we conclude that $%
x=H_{1}\ast H_{2}\ast h_{n_{2}+1}\ast ...\ast h_{r}$ where the rest of the $%
h_{i}$ are disjoint with $H_{1}\ast H_{2}$. Repeat the first step with $%
h_{n_{2}+1}$ and so on to get $x=H_{1}\ast H_{2}\ast H_{3}\ast ...\ast H_{n}$
where each $H_{i}$ disjoint with all the previous ones. Because similarity
is an equivalence relation $H_{i}$ are mutually disjoint, being separated on
the basis of similarity and disjointness and of course each of $H_{i}$ is
homogeneous by construction.

There's an alternative method that may perhaps be easier for some, though
harder in practice. Select from 
\begin{equation*}
H(x)=\{h_{1},...,h_{r}\}
\end{equation*}%
a set 
\begin{equation*}
\mathcal{K}_{k}=\{h_{01},h_{02},...,h_{0k}\}
\end{equation*}%
of mutually disjoint homogeneous factors/summands of $x$. If there is a
member of $H(x)$ that is disjoint with each member of $\mathcal{K}_{k}$,
then label it as $h_{0k+1}$ and form $\mathcal{K}_{k+1}=%
\{h_{01},h_{02},...,h_{0k+1}\}$. Repeat until you get to a stage 
\begin{equation*}
\mathcal{K}_{n}=\{h_{01},h_{02},...,h_{0n}\}
\end{equation*}%
where $h_{0i}$ are mutually disjoint and there's no member $h_{j}$ of $H(x)$
left such that $h_{j}\wedge h_{0i}=e$. In this case claim that we have a
maximal set of mutually disjoint homogeneous factors/summands of $x$. Now
suppose that there are two sets 
\begin{equation*}
\mathcal{K}_{m}=\{h_{01},h_{02},...,h_{0m}\}
\end{equation*}%
and 
\begin{equation*}
\mathcal{K}_{n}=\{k_{01},k_{02},...,k_{0n}\}
\end{equation*}%
where $m<n.$ As $\mathcal{K}_{m}$ is maximal, and as we are dealing with
homogeneous elements, each $k_{0i}$ is similar to exactly one of $h_{0j}$
and similars can replace similars. By relabeling we can assume that $k_{0i}$
replaces $h_{0i}.$ But then the extra ones in $\mathcal{K}_{n}$ would have
to be similar to some that they were disjoint with in $\mathcal{K}_{n}$.
Once that is settled take the maximal set $\mathcal{K}_{n}=%
\{h_{01},h_{02},...,h_{0n}\}$ and form: $H(h_{0i})=\{h\in H(x)|h\sim h_{0i}\}
$. Then $\{H(h_{0i})\}$ is a partition of $H(x)$. Next write $H_{i}=\Pi
_{\ast }h$ where $h$ varies over $H(h_{0i})$. Then by construction and by
Propositions \ref{Proposition YX5}, \ref{Proposition YX6} and \ref%
{Proposition YX9} $x=H_{1}\ast ...\ast H_{n}$ is a product of mutually
disjoint homogeneous elements. Now suppose that 
\begin{equation*}
x=H_{1}\ast ...\ast H_{n}=K_{1}\ast ...\ast K_{m}
\end{equation*}%
where $H_{1},...,H_{n}$ (resp., $K_{1},...,K_{m}$) are mutually disjoint
homogeneous elements. Since $H_{1}\leq x=K_{1}\ast ...\ast K_{m}$, $%
H_{1}\leq K_{j}$ for exactly one $j$, by Propositions \ref{Proposition YX8}.
Similarly as $K_{j}\leq x=H_{1}\ast ...\ast H_{n}$ and as, being
homogeneous, $K_{j}$ cannot be nondisjoint with two disjoint elements we
conclude that $K_{j}\leq H_{1}$. That is $K_{j}=H_{1}$. Repeating with $H_{2}
$ and so on we conclude that $n\leq m$. Similarly repeating the whole
process with $K_{i}$ we get $m\leq n$. Now relabeling $K_{j}$ we can get the
correspondence $H_{i}=K_{i}$.

If $M$ is a $\wedge $-smooth pre-Riesz monoid with a non-empty set $h(M)$ of
homogeneous elements we can always form a sub-semigroup 
\begin{equation*}
H^{\ast }(M)=\{h_{1}\ast h_{2}\ast ...\ast h_{n}|h_{i}\in h(M)\}
\end{equation*}%
of finite sums/products of members of $h(M).$ The sub-semigroup $H^{\ast }(M)
$ can be made into a monoid $H_{h}(M)=H^{\ast }(M)\cup \{e\}$ by throwing in 
$e$ as an empty sum/product. Call a directed p.o. monoid $M$ semi
homogeneous if $M$ is conic and each strictly positive element of $M$ is a
finite sum of homogeneous elements and call $M$ \emph{virtually factorial }(%
\emph{v-factorial}) if $M$ is semi-homogeneous such that every strictly
positive element of $M$ is uniquely expressible as a sum/product of finitely
many mutually disjoint elements.

\begin{corollary}
\label{Corollary YX11}Let $M$ be a $\wedge $-smooth pre-Riesz monoid with $%
h(M)$ $\neq \phi $. Then the submonoid 
\begin{equation*}
H_{h}(M)=\{h_{1}\ast h_{2}\ast ...\ast h_{n}|h_{i}\in h(M)\}\cup \{e\}
\end{equation*}
is a v-factorial monoid. Moreover $H_{h}(M)$ is a $\wedge $-smooth pre-Riesz
monoid. Consequently if a $\wedge $-smooth pre-Riesz monoid $M$ satisfies
CFC, then $H_{h}(M)$ is a v-factorial monoid.
\end{corollary}

\begin{proof}
That $H_{h}(M)$ is v-factorial follows from Theorem \ref{Theorem YX10}. For
the pre-Riesz part let $e<x,y\in H(M),$ under the induced partial order from 
$M,$ and suppose that $x\wedge y\neq e.$ Since $x\wedge y\neq e$ means there
is $g\leq x,y$ with $g\nleq e$ and since $H_{h}(M)$ is conic by construction
we conclude that there is $e<t\leq x,y$ in $H_{h}(M).$ The $\wedge $-smooth
part follows because the partial order is induced. The consequently part is
obvious.
\end{proof}

This brings us to the examples of virtual factoriality. We would be
selective, as they are strewn all over. For a start let us recall that a
nonzero finitely generated ideal $A$ of a domain $D$ is called primitive if $%
A\subseteq aD$ implies that $a$ is a unit, for each $a\in D$ and $A$ is
called super-primitive if $A_{v}=D$. Call $D$ a PSP (primitive is
super-primitive) domain if every primitive finitely generated ideal of $D$
is super-primitive. It was established in \cite{YZ 2011} that $D$ is a PSP
domain if and only if $G(D)$, the group of divisibility of $D$, is a
pre-Riesz group. Indeed, in light of what we have established in this paper, 
$D$ is a PSP domain if the monoid of nonzero principal ideals of $D$ is a
pre-Riesz monoid. Also, as a Riesz monoid is a pre-Riesz monoid and a
pre-Schreier domain is a domain $D$ whose monoid of nonzero principal ideals
is a Riesz monoid and it has been established, often, that a pre-Schreier
domain is a PSP domain, (see \cite{MR 1978}, \cite{Zaf 1987}, etc.) Finally
let's recall that $r\in D\backslash \{0\}$ is called \emph{rigid} if for all 
$x,y|r$ we have $x|y$ or $y|x$ and as in \cite{Zaf 2021} let's call $D$
semirigid if every nonzero non unit of $D$ is expressible as a finite
product of rigid elements.

\begin{example}
\label{Example ZX1} (1) Let $D$ be a semirigid PSP domain then $G(D)^{+}$ or 
$m(D)$ the monoid of nonzero principal ideals of $D$ is a v-factorial
monoid. Of course if $D$ is a GCD domain and semirigid we get a semirigid
GCD domain of \cite{Zaf 2022} where $m(D)$ is a factorial monoid. Example
3.7 of \cite{Zaf 2022} serves as an example of a semirigid Schreier domain.
Of course a UFD is a PSP domain and can be treated as a semirigid PSP domain.
\end{example}

\begin{example}
Let $D$ be a PSP domain of finite $t$-character, then the monoid $%
H_{h}(m(D)) $ generated by homogeneous elements of $m(D)=\{xD|x\in
D\backslash \{0\}\}$ is a v-factorial monoid. ($D$ is of finite $t$%
-character if every nonzero non unit of $D$ is contained in at most a finite
number of maximal $t$-ideals of $D$.$)$
\end{example}

\begin{example}
Let $D$ be a PSP domain such that every maximal $t$-ideal $M$ of $D$ is
associated to a homogeneous element $h$ of $D$, i.e $M=M(h)=\{x\in
D|(x,h)_{v}\neq D\}$. Then $H_{h}(m(D))$ is a v-factorial monoid.
\end{example}

Illustration: Let for $xD,yD\in m(D)$,$xD\wedge yD$ stand for $(xD+yD)_{v}$.
Then $xD,yD$ are disjoint if $(xD+yD)_{v}=D$ and non-disjoint if $%
(xD+yD)_{v}\neq D$. Obviously $m(D)$ is pre-Riesz if for every finite set $%
x_{1},...,x_{n}$, $(\sum x_{i}D)_{v}=D$ or $\sum x_{i}D\subseteq dD$ for
some nonzero non unit $d$ and for $x\in D\backslash \{0\}$,$xD$ is
homogeneous if for each pair of nonzero principal ideals containing $xD$ are
non-disjoint (i.e. for all non units $u,v|x$, $u,v$ have a non unit common
divisor). Of course a rigid element has this property. So if (1) $D$ is a
semirigid PSP (pre-Schreier or GCD) domain, all of $m(D)$ is v-factorial.

For (2) recall from section 4 of \cite{YZ 2011} that a PSP domain $D$ is of
finite $t$-character if and only if every nonzero non unit of $D$ is
divisible by at least one homogeneous element and by at most a finite number
of homogeneous elements. Once we note that $h(m(D))$, the set of homogeneous
elements of $m(D)$, is non-empty, we conclude, via Corollary \ref{Corollary
YX11}, that $H_{h}(m(D))$ is a v-factorial monoid. The reader can find
several examples of non-pre-Schreier PSP domains of finite $t$-character in
section 3, e.g. Example 3.7, of \cite{Zaf 2021}.

For (3) our arguments are the same as the ones for (2) An example of a PSP
domain is $D=Z+XR[[X]]$, where $Z$ is the ring of integers and $R$ is the
field of real numbers, as explained in the illustration of Example 4.9 of 
\cite{YZ 2011}.

Let, for a star operation $\star $ of finite character, $\varphi _{\star }(D)
$ denote the set of (nonzero) $\star $-ideals of finite type. Then $\varphi
_{\star }(D)$, is a monoid under the product called $\star $-multiplication $%
\times ^{\star }:I\times ^{\star }J=(IJ)^{\star }$, where $I,J\in \varphi
_{\star }(D)$, while $\varphi _{\star }(D)$ is closed under $\star $-sum $%
+^{\star }:$ $I+^{\star }J=(I+J)^{\star }$, allowing for

\begin{center}
$K\times ^{\star }(I+^{\star }J)=K\times ^{\star }I+^{\star }K\times ^{\star
}J=(KI+KJ)^{\star }$.

Taking $I\leq J\Leftrightarrow I\supseteq J$ in $\varphi _{\star }(D)$, we
have $(I+^{\star }J)=\inf (I,J)=I\wedge J$. Consequently, we have a $\wedge $%
-semilattice because $+^{\star }$ is indeed associative. With the order
defined as we have $D$ is the least element of $\varphi _{\star }(D)$. That $%
\varphi _{\star }(D)$ is $\wedge $-smooth because generally in $\varphi
_{\star }(D)$ we have $K\times ^{\star }I\wedge K\times ^{\star }J)=K\times
^{\star }(I\wedge J)$.

Now requiring that whenever, for $I_{1},...,I_{n}\in \varphi _{\star }(D)$,
the proper ideal $(I_{1},...,I_{n})^{\star }$ exceeds some member of $%
\varphi _{\star }(D)$, which indeed is the case, we make $\varphi _{\star
}(D)$ into a $\wedge $-smooth pre-Riesz monoid, with $e=D$.
\end{center}

The idea took root in \cite{DZ 2010}, when pre-Riesz monoids were not in
plain sight, with just the notion of factoriality in Riesz groups of \cite%
{MRZ 2008} and some of the work on factorization to go on. A nonzero $\star $%
-ideal $\mathbf{h}$ of finite type was called, in \cite{DZ 2010}, a
homogeneous ideal if for every pair $\mathbf{J,k}$ of proper $\star $-ideals
of finite type containing $\mathbf{h}$ we had $(\mathbf{j,k})_{v}\neq D$ and
now, in more general terms, it is: $\mathbf{j\wedge k}\neq e$ for all $%
\mathbf{j,k}\in (e,\mathbf{h}]$. It was shown in \cite{DZ 2010} that $%
\mathbf{h\in }$ $\varphi _{\star }(D)$ is homogeneous if and only if $%
\mathbf{h}$ is contained in a unique maximal $\star $-ideal of $D$. Later,
in \cite{AZ 2019}, a domain $D$ was called $\star $-\emph{semi homogeneous} (%
$\star $-\emph{SH}) \emph{domain} if every nonzero principal ideal of $D$
was expressible as a $\star $-product of homogeneous ideals. According to
Theorem 6 of \cite{AZ 2019}: Given that $\star $ is a star operation of
finite character, if $I$ is a nonzero $\ast $-ideal of finite type in a $%
\star $-SH domain $D$ such that $I\neq D$, then $I$ is uniquely expressible
as a $\star $-product of mutually $\star $-comaximal $\star $-homogeneous
ideals. With the above description and illustration, the following results
and examples can be established. Indeed it was shown in \cite{DZ 2010} that $%
\varphi _{\star }(D)$ satisfied CFC if and only if $D$ was of finite $\star $%
-character, (see also \cite{Zaf 2021}).

\begin{proposition}
\label{Proposition ZX2}Let $\star $ be a finite character star operation
defined on an integral domain $D$. Then the monoid $\varphi _{\star }(D)$ of
nonzero $\star $-ideals of finite type of $D$, under $\star $%
-multiplication, is a $\wedge $-smooth pre-Riesz monoid with order defined
by $I\leq J$ if and only if $I\supseteq J$, for $I,J\in \varphi _{\star }(D)$%
. Moreover the following hold: (1) $D$ is a $\star $-SH domain if and only
if $\varphi _{\star }(D)$ is a v-factorial monoid, (2) If $h(\varphi _{\star
}(D))\neq \phi $, then $H_{h}(\varphi _{\star }(D))$ is a v-factorial monoid.
\end{proposition}

\textbf{Acknowledgements:} I am indebted to Brian Davey, the managing editor
of Algebra Universalis, for straightening my original write up to make it
look like a decent paper.


\begin{thebibliography}{99}
\bibitem{And 1978} Anderson, D.D: $\pi $-domains, overrings and divisorial
ideals, Glasgow Math J. , \textbf{19}(2), 199-203 (1978)

\bibitem{ACZ 2015} Anderson, D.D, Chang, G.W., Zafrullah, M.: Nagata-like
theorems for integral domains of finite character and finite $t$-character.
J. Algebra Appl.\textbf{\ 14}(8) (2015)

\bibitem{ADZ 2007} Anderson, D.D, Dumitrescu, T., Zafrullah, M.:
Quasi-Schreier domains II, Comm. Algebra \textbf{35}, 2096-2104 (2007)

\bibitem{AZ 2019} Anderson, D.D., Zafrullah, M.: On $\star $-semi
homogeneous integral domains, in Advances in Commutative Algebra, Editors:
A. Badawi and J. Coykendall, pp. 7-31 (2019)

\bibitem{B 1948} Birkhoff, G.: Lattice Theory, Amer. Math. Soc. Colloq.
Publication \textbf{25}(1948)

\bibitem{Coh 1968} Cohn, P.M.: Bezout rings and their subrings, Proc.
Cambridge Philos. Soc.\textbf{\ 64}, 251-264 (1968)

\bibitem{Con 1961} Conrad, P.: Some structure theorems for lattice ordered
groups, Trans. Amer. Math. Soc. \textbf{99}, 212-240 (1961)

\bibitem{CMZ 1978} Costa, D.L., Mott, J.L., Zafrullah, M.: The construction $%
D+XD_{S}[X]$,\ J. Algebra \textbf{53}, 423-439 (1978)

\bibitem{CMZ 1986} Costa, D.L., Mott, J.L., Zafrullah, M.: Overrings and
dimensions of general $D+M$ constructions,\ J. Natur. Sci. and Math. \textbf{%
26}(2), 7-14 (1986)

\bibitem{DM 2003} Dumitrescu, T., Moldovan, R.: Quasi-Schreier domains,
Math. Reports \textbf{5}, 121-126 (2003)

\bibitem{DZ 2010} Dumitrescu, T., Zafrullah, M.: Characterizing domains of
finite $\star $-character, J. Pure Appl. Algebra, \textbf{214}, 2087-2091
(2010)

\bibitem{DZ 2011} Dumitrescu, T., Zafrullah, M.: $t$-Schreier domains, Comm.
Algebra, \textbf{39}, 808-818 (2011)

\bibitem{EIT 2019} El-Baghdadi, S., Izelgue, L., Tamoussit, A.: Almost Krull
domains and their rings of integer-valued polynomials, J. Pure Appl.
Algebra, 224(6) 106269, 9 pp (2020)

\bibitem{Fos 1973} Fossum, R.: The divisor class group of a Krull domain,
Ergebnisse der Mathematik und ihrer grenzgebiete B. \textbf{74},
Springer-Verlag, Berlin, Heidelberg, New York, (1973)

\bibitem{Fu 1965} Fuchs, L.: Riesz groups, Ann. Sc. Norm. Super. Pisa, Sci.
Fis. Mat., III. Ser. \textbf{19}, 1--34 (1965)

\bibitem{GBL 2004} Gabelli, S., Houston, E., Lucas, T.: The t\#-property for
integral domains, J. Pure Appl. Algebra \textbf{194}, 281-298 (2004)

\bibitem{Gil 1972} Gilmer, R.: Multiplicative Ideal Theory, Marcel-Dekker
(1972)

\bibitem{Gil 1984} Gilmer, R.: Commutative Semigroup Rings, Chicago Lectures
in Mathematics, University of Chicago Press, Chicago, IL, 1984

\bibitem{GV 1977} Glaz. S. , Vascocelos, W.: Flat ideals II,\ Manuscripta
Math. \textbf{22}, 325-341 (1977)

\bibitem{FHK B} Halter-Koch, F.: Ideal Systems, An introduction to ideal
theory, Marcel Dekker, New York, (1998)

\bibitem{FHK mix} Halter-Koch, F.: Mixed invertibility and Prufer-like
monoids and domains, \textquotedblright\ in 
\"{}
Commutative Algebra and Its Applications, Proceedings of the Fez Conference,
pp. 247--258, Walter de Gruyter, Berlin, Germany, (2009)

\bibitem{HO 1973} Heinzer, W., Ohm, J.: An essential ring which is not a
v-multiplication ring, Canad. J. Math.\textbf{\ 25}, 856-861 (1973)

\bibitem{HZ 1988} Houston, E., Zafrullah, M.: Integral domains in which
every $t$-ideal is divisorial, \ Michigan Math. J. \textbf{35}, 291-300
(1988)

\bibitem{HZ 1989} Houston, E., Zafrullah, M.: On $t$-invertibility II,\
Comm. Algebra \textbf{17}(8), 1955-1969 (1989)

\bibitem{HZ 2015} Houston, E., Zafrullah, M.: Integral domains in which any
two $v$-coprime elements are comaximal, J. Algebra \textbf{423}, 93-115
(2015)

\bibitem{Kap 1970} Kaplansky, I.: Commutative Rings, Allyn and Bacon, Boston
(1970).

\bibitem{LYZ 2020} Liu, Pengyuan , Yang, Y.C., Zafrullah, M.: Conrad's
F-condition for partially ordered monoids, Soft Comput \textbf{24},
9375--9381 (2020)

\bibitem{MR 1978} McAdam, S., Rush, D.E.: Schreier rings, Bull. London Math.
Soc, \textbf{10}, 77-80 (1978)

\bibitem{Mat 1988} Matsuda, R.: Notes on torsion-free Abelian semigroup
rings, Bull. Fac. Sci., Ibaraki Univ., \textbf{20}, 51-59 (1988)

\bibitem{Mat 1997} Matsuda, R.: Note on Schreier semigroup rings, Math. J.
Okayama Univ. \textbf{39}, 41--44 (1997)

\bibitem{MRZ 2008} Mott, J., Rashid, M., Zafrullah, M.: Factoriality in
Riesz groups, J. Group Theory \textbf{11} (1), 23-41 (2008)

\bibitem{MZ 1981} Mott, J., Zafrullah, M.: On Prufer v-multiplication
domains, Manuscripta Math. \textbf{35}, 1-26 (1981)

\bibitem{MZ 1990} Mott, J., Zafrullah, M.: Unruly Hilbert domains, Canad.
Math. Bull. \textbf{33} (1), 106-109 (1990)

\bibitem{MZ 1991} Mott, J., Zafrullah, M.: On Krull domains, Arch. Math. 
\textbf{56}, 559-568 (1991)

\bibitem{Nishi 1973} Nishimura, T: On $t$-ideals of an integral domain, J.
Math. Kyoto Univ. (JMKYAZ) \textbf{13}-1, 59-65 (1973)

\bibitem{YZ 2011} Yang, Y.C., Zafrullah, M.: Bases of pre-Riesz groups and
Conrad's F-condition, Arab. J. Sci. Eng. \textbf{36}, 1047--1061 (2011)

\bibitem{Zaf 1986} Zafrullah, M.: On generalized Dedekind domains,
Mathematika \textbf{33}, 285-295 (1986)

\bibitem{Zaf 1987} Zafrullah, M: On a property of pre-Schreier domains,
Comm. Algebra \textbf{15}, 1895-1920 (1987)

\bibitem{Zaf 2000} Zafrullah, M.: Putting $t$-invertibility to use,
Non-Noetherian Commutative Ring Theory, in: Math. Appl., vol. \textbf{520},
Kluwer Acad. Publ., Dordrecht, pp. 429-457 (2000)

\bibitem{Zaf 2006} Zafrullah, M.: What $v$-coprimality can do for you, in
Multiplicative ideal theory in commutative algebra, 387--404, Springer, New
York, (2006)

\bibitem{Zaf 2020} Zafrullah, M.: Domains whose ideals meet a universal
restriction, https://arxiv.org/pdf/2006.04135.pdf

\bibitem{Zaf 2021} Zafrullah, M.: On $\star $-potent domains and $\star $%
-homogeneous ideals (accepted for publication in: Rings, Monoids, and Module
Theory (Eds. A. Badawi and J. Coykendall), Springer (to appear).. Also
available at : https://arxiv.org/pdf/1907.04384.pdf

\bibitem{Zaf 2022} Zafrullah, M.: Semirigid GCD domains II, J. Algebra and
Applications, (to appear) DOI: 10.1142/s0219498822501614 (2022)
\end{thebibliography}
\end{document}